\documentclass[11pt]{article}
\usepackage{graphicx}
\usepackage{amsmath,amssymb,amsthm,amsfonts}
\usepackage{color}
\usepackage[utf8]{inputenc}  
\usepackage[T1]{fontenc}     
\usepackage{amssymb}
\usepackage{cite}
\newtheorem{thm}{Theorem}[section]
\newtheorem{cor}[thm]{Corollary}
\newtheorem{lem}[thm]{Lemma}
\newtheorem{prop}[thm]{Proposition}

\theoremstyle{definition}
\newtheorem{defn}{Definition}[section]
\theoremstyle{remark}
\usepackage{appendix}
\newtheorem{rem}{Remark}[section]

\numberwithin{equation}{section}

\DeclareMathSymbol{\C}{\mathalpha}{AMSb}{"43}

\newcommand{\eps}{\varepsilon}

\newcommand{\R}{{\mathbb{R}}}

\newcommand{\inte}{\int_{\mathbb{R}}}

\def\R{{\mathbb R}}
\def\C{{\mathbb C}}
\allowdisplaybreaks[4]

\newcommand{\bsub}{\begin{subequations}}
\newcommand{\esub}{\end{subequations}$\!$}

\begin{document}
	
%
%
%
%
%
%
%
%

\title{Classification and Nondegeneracy of Cubic Nonlinear Schr\"{o}dinger System in $\R$}

\author{Yujin Guo\thanks{School of Mathematics and Statistics,  Key Laboratory of Nonlinear Analysis $\&$ Applications (Ministry of Education), Central China Normal University, Wuhan 430079, P. R. China. Y. J. Guo is partially
supported by  NSF of China
(Grants 12225106 and 12371113) and National Key R $\&$ D Program of China (Grant 2023YFA1010001). Email: yguo@ccnu.edu.cn. }, \  Yong Luo\thanks{School of Mathematics and Statistics, and Hubei Key Laboratory of
Mathematical Sciences,  Central China Normal University, Wuhan 430079, P. R. China. Email: yluo@ccnu.edu.cn.},\ \,  and\; Juncheng Wei\thanks{Department of Mathematics, Chinese University of Hong Kong, Shatin, NT, Hong Kong. J. Wei is partially supported  by   Hong Kong General Research Fund ``New frontiers in singular limits of nonlinear partial differential equations". Email: wei@math.cuhk.edu.hk.}  }

\date{}

\smallbreak \maketitle
\begin{abstract}
We study the following one-dimensional cubic nonlinear Schr\"{o}dinger system:
\[
u_i''+2\Big(\sum_{k=1}^Nu_k^2\Big)u_i=-\mu_iu_i \ \,\ \mbox{in}\, \ \mathbb{R} , \ \  i=1, 2, \cdots, N,
\]
where $\mu_1\leq\mu_2\leq\cdots\leq\mu_N<0$ and $N\ge 2$. In this paper, we mainly focus on the case $N=3$ and prove the following results: (i). The solutions of the system can be completely classified; (ii). Depending on the explicit values of $\mu_1\leq\mu_2\leq\mu_3<0$, there exist two different classes of normalized solutions $u=(u_1, u_2, u_3)$ satisfying $\int _{\R}u_i^2dx=1$ for all $i=1, 2, 3$, which are completely different from the case $N=2$; (iii). The linearized operator at any nontrivial solution of the system is non-degenerate. The conjectures on the explicit classification and nondegeneracy of solutions for the system are also given for the case $N>3$. These address the questions of [R. Frank, D. Gontier and M. Lewin, CMP, 2021], where the complete classification and uniqueness results for the system were already proved for the case $N=2$.
\end{abstract}

\vskip 0.05truein


\vskip 0.2truein

\section{Introduction}

In this paper, we consider the following one-dimensional cubic  nonlinear Schr\"{o}dinger system
\begin{equation}\label{GPN}
		u_i''+2\Big(\sum_{k=1}^Nu_k^2\Big)u_i=-\mu_iu_i \ \,\ \mbox{in}\, \ \R , \ u_i \in H^1(\R) \ \  i=1, 2, \cdots, N,
\end{equation}
where $\mu_1\leq\mu_2\leq\cdots\leq\mu_N<0$ and $N\ge 2$.
The nonlinear Schr\"{o}dinger system (\ref{GPN}) arises from various physical situations (cf. \cite{i,M}), such as multiple-component Bose-Einstein condensates (cf. \cite{LS,LW}), nonlinear optics (cf. \cite{ZS}), Langmuir waves  (cf. \cite{FI}), and so on. Especially, the system (\ref{GPN}) can be used to model the ground states of $N$-component Bose gases, see \cite{LS,LW} and the references therein. On the other hand, the system (\ref{GPN}) can be also employed to discuss the  ground states of $N$-particle fermionic quantum systems, which were recently investigated in \cite{Lewin,i,Lewin3} and the references therein.

Another motivation of studying (\ref{GPN}) is the relation with Lieb-Thirring inequality.  \textcolor{black}{Indeed, if $\mu_i\neq \mu_j$, then $u_i$ and $u_j$ are eigenfunctions  of the operator $$H:=-\frac{d^2}{dx^2}-2\sum_{n=1}^Nu_n^2 \ \,\ \mbox{in}\ \, L^2(\R),$$
and hence they are automatically orthogonal.}
We further note from \cite{Lewin,Lewin3} and Remark \ref{rem-lt} below that the solutions $u=(u_1, u_2, \cdots, u_N)$ of the system (\ref{GPN}) are connected  with the following finite-rank Lieb-Thirring inequality
\begin{equation}\label{lt}
		 \sum_{n=1}^N\big|\mu_n\big(-\Delta -V(x)\big)\big|^\gamma\le L^{(N)}_{\gamma, d} \int_{\R^d}V(x)_+^{\gamma +\frac{d}{2}}dx , \ \  V(x)\in L^{\gamma +\frac{d}{2}}(\R^d),
\end{equation}
in the critical case $\gamma=\frac{1}{2}$ and $d=1$, where $V(x)_+=\max\{V(x), 0\}$ and $\mu_n\big(-\Delta -V(x)\big)$ denotes the $nth$ min-max level of $-\Delta -V(x)$,
so that it equals to the $n$th negative eigenvalue (counted with multiplicity) when it exists, and $0$ otherwise. By the definition of (\ref{lt}), one can observe   that the best constant $L^{(N)}_{\gamma, d}$ of (\ref{lt}) satisfies $L^{(N)}_{\gamma, d}\le L^{(N+1)}_{\gamma, d}$   for all $N\in \mathbb{N}^+$. Especially, it was proved in \cite{W,HLT} that  for the critical case $\gamma=\frac{1}{2}$ and $d=1$, the best constant $L^{(N)}_{\gamma, d}$ satisfies
\[
L_{\frac{1}{2}, 1}:=L^{(\infty)}_{\frac{1}{2}, 1}=\lim_{N\to\infty}L^{(N)}_{\frac{1}{2}, 1}=\frac{1}{2},
\]
where the constant $L_{\gamma, d}$ is called the best Lieb-Thirring constant, see \cite{W,HLT,Lewin,i,Lewin3} and the references therein.

When $N=1$ and $\mu_1=-1$, direct integration yields that up to the translation and the sign $\pm,$  \eqref{GPN} admits a unique solution
\begin{equation}\label{2.0M}
u(x)\equiv \frac{1}{\cosh(x)}.
\end{equation}
When $N=2$, the system \eqref{GPN} was analyzed recently in  \cite{Lewin}, where the authors obtained the complete classification on the general solutions of \eqref{GPN} by applying Hirota's bilinearisation method (cf. \cite{RL,RSL,H}). Based on the classification on the general solutions of \eqref{GPN}, the authors in \cite{Lewin} derived successfully the following interesting uniqueness, up to translations and the uncorrelated signs $\pm$, of normalized solutions for \eqref{GPN} in the case $N=2$.

\vskip 0.05truein

\noindent{\bf Theorem A}{\em \ (\cite[Theorem 8]{Lewin})  For $\mu_1\leq\mu_2<0$, suppose $u=(u_1,u_2)\in H^1(\R)^2$ is a normalized solution of \eqref{GPN} with $N=2$, in the sense that $\|u_1\|_{L^2(\R)}=\|u_2\|_{L^2(\R)}=1$. Then it necessarily has $\mu_1=\mu_2<0$, and up to translations, $u$ satisfies
\begin{equation}\label{2.1}
	u_1(x)\equiv\pm\frac{\sqrt{2}}{2\cosh(x)},\,\, u_2(x)\equiv\pm\frac{\sqrt{2}}{2\cosh(x)}
\end{equation}
for two uncorrelated signs $\pm$.	}


The authors also conjectured in \cite[Subsection 1.3]{Lewin} that the above mentioned results can be probably generalized to the system \eqref{GPN} for any $N\ge 3$. On the other hand, the existing results show that it is usually very useful to know the complete classification and  nondegeneracy of solutions for the system \eqref{GPN}, which unfortunately seem very challenging for general situations, see \cite{A,BS,LW,LWY} and the references therein. Stimulated by these facts, the main purpose of this paper is to address the complete classification, the nondegeneracy and some other analytical properties of solutions for the system \eqref{GPN} with $N\ge 3$. It turns out that the analysis of the system \eqref{GPN} with $N\ge 3$ is more complicated, for which one requires  new approaches. More precisely, in this paper we are able to understand completely the case $N=3$ of the system (\ref{GPN}), and we finally discuss some partial results, together with several conjectures, for the general case $N>3$ of the system (\ref{GPN}).


In this paper, we first focus on the case $N=3$ of the system (\ref{GPN}), which can be rewritten explicitly as the following form: 	
\begin{equation}\label{GPs}
	\left\{\begin{array}{lll}
		u_1''+2(u_1^{2}+u_2^{2}+u_3^{2})u_1=-\mu_1u_1 \ \ \mbox{in} \ \, \R,\\[3mm]
	u_2''+2(u_1^{2}+u_2^{2}+u_3^{2})u_2=-\mu_2u_2 \ \ \mbox{in} \ \, \R,\\[3mm]
		u_3''+2(u_1^{2}+u_2^{2}+u_3^{2})u_3=-\mu_3u_3 \ \ \mbox{in} \ \, \R,
	\end{array}\right.
\end{equation}
where $\mu_1\leq\mu_2\leq\mu_3<0$ are arbitrary. Based on Hirota's bilinearisation method (cf. \cite{H,RL,RSL,Lewin}), we shall employ  the algebraic analysis to prove the following classification on general solutions of \eqref{GPs}.

\begin{thm}\label{thm1.1}
For $\mu_1\leq\mu_2\leq\mu_3<0$, suppose $u=(u_1,u_2,u_3)$ is a  solution of \eqref{GPs} in $H^1(\R)^3$. Then there exists a vector $(a_1,a_2,a_3)\in\R^3$ such that $u=(u_1,u_2,u_3)$ can be written exactly as:
\begin{equation}\label{exp1}
u_i(x)=\frac{g_{i}(x)}{f(x)},\ \ i=1,2,3,
\end{equation}
where $g_i(x)$ and $f(x)$ satisfy for $\eta_i=\sqrt{|\mu_i|}>0$,
\[
\begin{aligned}
g_i(x)&=a_ie^{\eta_ix}+\sum_{j=1}^{3}\frac{a_ia_j^2(\eta_i-\eta_j)}{4\eta_j^2(\eta_i+\eta_j)}e^{(2\eta_j+\eta_i)x}\\
&\quad+\sum_{1\leq j<k\leq 3}\frac{a_ia_j^2a_k^2(\eta_j-\eta_k)^2(\eta_i-\eta_j)(\eta_i-\eta_k)}
{16\eta_j^2\eta_k^2(\eta_j+\eta_k)^2(\eta_i+\eta_j)(\eta_i+\eta_k)}e^{(\eta_i+2\eta_j+2\eta_k)x},\ \ i=1,2,3,
\end{aligned}
\]
 and
\[
\begin{aligned}
f(x)&=1+\sum_{j=1}^{3}\frac{a_j^2}{4\eta^2_j}e^{2\eta_jx}+\sum_{1\leq j<k\leq 3}\frac{a_j^2a_k^2(\eta_j-\eta_k)^2}{16\eta_j^2\eta_k^2(\eta_j+\eta_k)^2}e^{2(\eta_j+\eta_k)x}\\
&\quad+\frac{a_1^2a_2^2a_3^2(\eta_1-\eta_2)^2(\eta_1-\eta_3)^2(\eta_2-\eta_3)^2}{64\eta_1^2\eta_2^2\eta_3^2(\eta_1+\eta_2)^2(\eta_1+\eta_3)^2(\eta_2+\eta_3)^2}e^{2(\eta_1+\eta_2+\eta_3)x}.
\end{aligned}
\]
\end{thm}

\textcolor{black}{
\begin{rem} \label{rem-O}
In Figure 1 below, we plot the function $V(x)=2(u_1^2+u_2^2+u_3^2)$  with fixed parameters $(\mu_1,\mu_2,\mu_3)=(-4,-2,-1)$ and $(a_2,a_3)=(1,0.00674)$, where $a_1>0$ varies. It shows in Figure 1 that the maximal point of $V(x)$  moves to the left as $\ln a_1$ increases, and $V(x)$  exhibits three bumps. It seems from Figure 1 that the solutions of Theorem \ref{thm1.1} have the properties similar to the one of ``solitons''.
  \begin{figure}[t]
		\centering
		\includegraphics[width=0.9\textwidth]{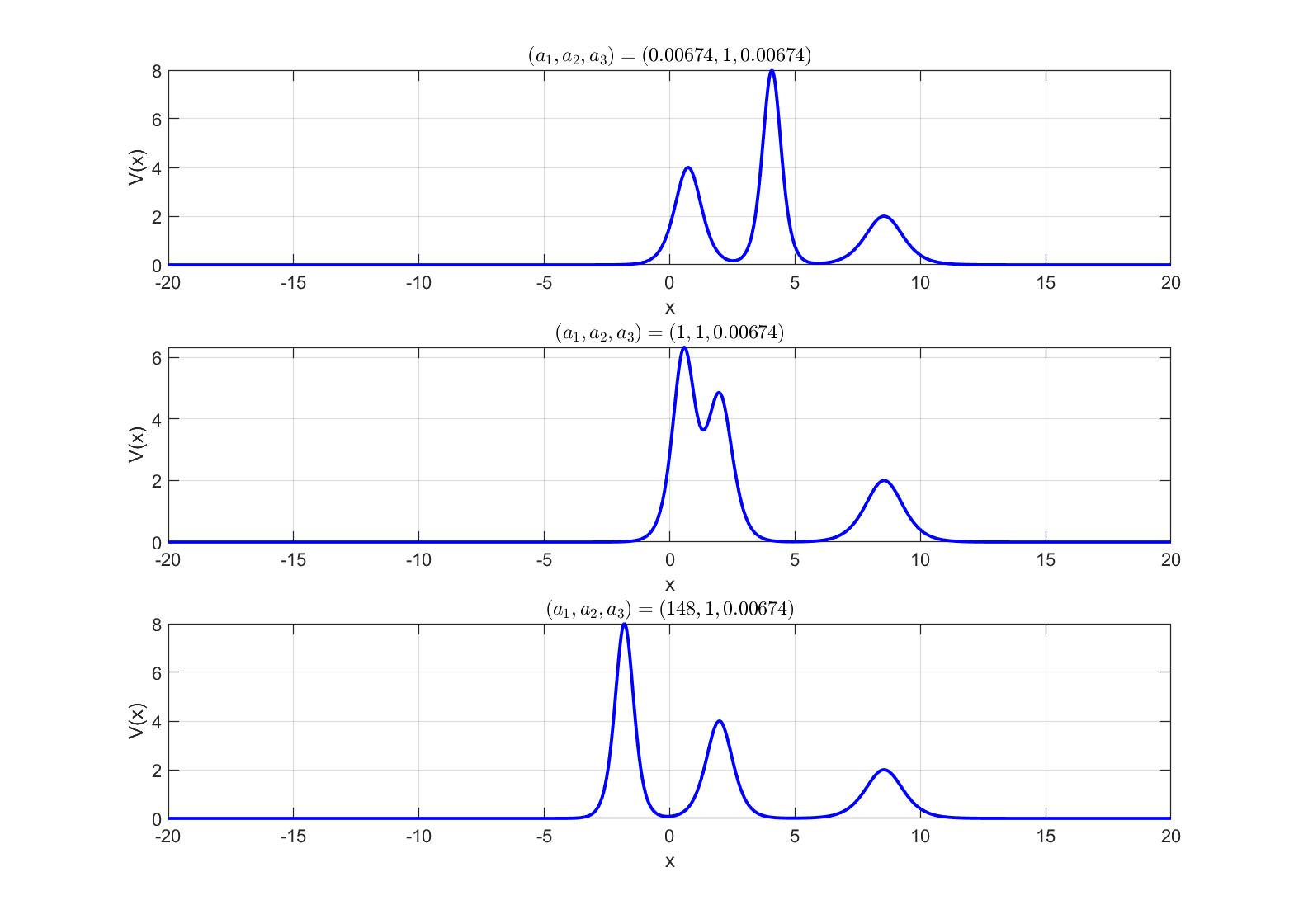}
		\caption{The curves of $V(x)=2(u_1^2+u_2^2+u_3^2)$ with fixed parameters $(\mu_1,\mu_2,\mu_3)=(-4,-2,-1)$ and $(a_2,a_3)=(1,0.00674)$, where $a_1>0$ varies.}
		{}
		\label{fig-O1} %
	\end{figure}
\end{rem}  }

\begin{rem} \label{rem-lt}
In order to prove Theorem \ref{thm1.1}, we shall derive in Lemma \ref{lem-2} that if $\mu_1<\mu_2<\mu_3<0$, then the solution $(u_1,u_2,u_3)\in H^1(\R)^3$ of \eqref{GPN} with $N=3$ satisfies $$\int_{\R}u^2_n(x)dx=2\eta_n=2\sqrt{|\mu_n|}, \ \ n=1,2,3.$$
This further implies that  the solution $(u_1,u_2,u_3)\in H^1(\R)^3$ of \eqref{GPN} with $N=3$ satisfies
\begin{equation}\label{rem-1:1}
	 \sum_{n=1}^3 |\mu_n |^\frac{1}{2}=\frac{1}{4} \int_{\R}V(x)dx, \ \ V(x)=2\sum_{n=1}^3u_n^2,
\end{equation}
which seems related to the finite-rank Lieb-Thirring inequality (\ref{lt}).
\end{rem}

\textcolor{black}{
\begin{rem} \label{rem-2}
It yields from Theorem \ref{thm1.1} that if $\mu_1\leq \mu_2\leq \mu_3<0$, then  $|u_1|>0$ holds in $\R$. This further implies that for $\mu_1\leq \mu_2\leq \mu_3<0$, $\mu_1$ is the first eigenvalue of the operator $H=-\frac{d^2}{dx^2} -2\sum_{n=1}^3u_n^2\ $ in $L^2(\R)$, and hence $\mu_1$ is simple by Krein-Rutman theorem. Moreover,
by the proof of Theorem \ref{thm1.1}, we can deduce from \eqref{sp-12:3} and \eqref{sp-13:3} below that  for $\mu_1\leq \mu_2\leq \mu_3<0$,
\[
\inte (u_1^2+u_2^2+u_3^2)^2dx=-\frac{2}{3}(\mu_1+\mu_2+\mu_3).
	\]	
Applying  the   finite-rank Lieb-Thirring inequality \eqref{lt} with $\gamma=\frac 3 2$ and $d=1$, one can further deduce from the above identity that if $\mu_1\leq \mu_2\leq \mu_3<0$ and $|\mu_1|>0$ is sufficiently larger than $|\mu_2|$ and $|\mu_3|$, then $\mu_2$ and $\mu_3$ are not the second and the third negative eigenvalues (counted with multiplicity) of the operator
$H=-\frac{d^2}{dx^2} -2\sum_{n=1}^3u_n^2\ $ in $L^2(\R)$.
\end{rem}
}

Theorem  \ref{thm1.1} gives a similar classification  of \cite[Lemma 15]{Lewin}, but the proof of Theorem  \ref{thm1.1} needs the more involved algebraic  analysis. Actually, applying Hirota's bilinearisation method (cf. \cite{H,RL,RSL}), we shall first prove that any function $u=(u_1,u_2,u_3)$ defined by \eqref{exp1} is a solution of \eqref{GPs}. We then derive  the associated integral system (\ref{I-1})--(\ref{I-3}) by applying the symmetry of \eqref{GPs}. Following (\ref{I-1})--(\ref{I-3}), if $\mu_2= \mu_i$ holds for either $i=1$ or $i=3$ or $i=1, 3$, see Proposition \ref{prop2A},  then the proof of Theorem  \ref{thm1.1} follows essentially from the argument of \cite[Lemma 15]{Lewin}. However, this unfortunately  is not applied to the challenging general case $\mu_1<\mu_2<\mu_3<0$, see Subsection 2.1 for more details.
Therefore, by algebraically analyzing the integral curves of the system \eqref{GPs}, we shall complete in Section 3 the proof of Theorem  \ref{thm1.1} for the general case where $\mu_1<\mu_2<\mu_3<0$.

In the following we discuss several applications of Theorem \ref{thm1.1}. Since the normalized solutions of the  system \eqref{GPs} were studied widely in recent few years, see \cite{A,BS} and the references therein, we shall first employ Theorem \ref{thm1.1} to establish the following refined classification on the normalized solutions of \eqref{GPs}:

\begin{thm}\label{thm1.2}
Suppose $\mu_1\leq\mu_2\leq\mu_3<0$. Then \eqref{GPs} admits normalized solutions $u=(u_1,u_2,u_3)\in H^1(\R)^3$ satisfying  $\int_{\R}u_i^2dx=1$ for all $i=1,2,3$, if and only if one of the following conditions holds:
\begin{itemize}
\item[(i).] $\mu_1=\mu_2=\mu_3=-\frac{9}{4}$;
\item[(ii).]    $\mu_1=\mu_2=-1$ and $\mu_3=-\frac{1}{4}$.
\end{itemize}
Moreover, if (i) holds, then  up to translations, \eqref{GPs} admits a unique normalized solution, in the sense that	
\begin{equation}\label{2.0}
	u_1(x)\equiv\pm\frac{\sqrt{3}}{2\cosh(\frac{3x}{2})},\,\, u_2(x)\equiv\pm\frac{\sqrt{3}}{2\cosh(\frac{3x}{2})},\,\,u_3(x)\equiv\pm\frac{\sqrt{3}}{2\cosh(\frac{3x}{2})}
\end{equation}
hold for three uncorrelated signs $\pm$. If (ii) holds, then up to translations, \eqref{GPs} admits infinitely many normalized solutions, which satisfy
\begin{equation}\label{thm1:M8}
\left\{\begin{split}
u_1(x)&=\frac{A e^{x}}{f(x)}\Big(1+\frac{B^2}{3}e^x\Big),\\
u_2(x)&=\pm  \frac{Ae^{x}}{f(x)}\Big(1+\frac{B^2}{3}e^x\Big),\\
u_3(x)&=\frac{B e^{\frac{1}{2}x}}{f(x)}\Big(1-\frac{A^2}{6}e^{2x}\Big),
\end{split}\right.
\end{equation}
and
\[
f(x)=1+B^2e^{x}+\frac{A^2}{2}e^{2x}+\frac{A^2B^2}{18}e^{3x},\,\ \forall A\neq 0,\forall B\neq 0.
\]
\end{thm}
Theorem \ref{thm1.2} shows that one cannot expect generally the uniqueness  of normalized solutions for the system \eqref{GPs}, which depends on the exact value of $\mu_1\leq\mu_2\leq\mu_3<0$.

\textcolor{black}{
\begin{defn}\label{Def1} The vector function
$u=(u_1,u_2,u_3)\in H^1(\R)^3$ is called a least energy normalized solution of  (\ref{GPs}), if $u=(u_1,u_2,u_3)\in H^1(\R)^3$ satisfies
\[
\begin{aligned}
E(u_1,u_2,u_3)=\inf\Big\{E(v_1,v_2,v_3):\ & (v_1,v_2,v_3)\,\  \hbox{satisfies (\ref{GPs}) for some $\mu_1\leq \mu_2\leq \mu_3<0$,}\\
&\hbox{and $\inte |v_i|^2dx=1$ for $i=1,2,3$}  \Big\},
\end{aligned}
\]
where the energy functional $E(v_1,v_2,v_3)$ is defined by
\[E(u_1,u_2,u_3):=	\sum_{k=1}^3\int_{\R}(u_k')^2dx-\int_{\R}\Big(\sum_{k=1}^3u_k^2\Big)^2dx.\]
\end{defn}
}

As other applications of Theorem \ref{thm1.1}, we further establish the following uniqueness of least energy normalized solutions  for (\ref{GPs}), and as well the nonexistence of orthonormal solutions $u=(u_1,u_2,u_3)\in H^1(\R)^3$ for (\ref{GPs}), in the sense that $(u_i, u_j)=\delta_{ij}$ holds for all $1\leq i\leq j\le 3$.

\begin{cor}\label{thm1.3}
		For $\mu_1\leq\mu_2\leq\mu_3<0$, we have the following conclusions:
\begin{itemize}		
\item[(i).] If $u=(u_1,u_2,u_3)\in H^1(\R)^3$ is a least energy normalized solution of \eqref{GPs}, then it necessarily has $\mu_1=\mu_2=\mu_3=-\frac{9}{4}$, and $u$ must satisfy \eqref{2.0}.

\item[(ii).] The system \eqref{GPs} does not  admit any orthonormal  solution $u=(u_1,u_2,u_3)\in H^1(\R)^3$.
\end{itemize}
\end{cor}

%

One can note from Theorem \ref{thm1.2} and Corollary \ref{thm1.3} that the least energy normalized solutions of \eqref{GPs} must be unique,  up to translations and the uncorrelated signs $\pm$. However, the number of higher energy normalized solutions for \eqref{GPs} must be either zero or  infinite, which depends on the exact value of $\mu_1\leq\mu_2\leq\mu_3<0$.

For any solution $u=(u_1,u_2,u_3)$ of the system \eqref{GPs}, we next consider the following linearized system of \eqref{GPs} at the solution $u$:
\begin{equation}\label{Lines}
	\left\{\begin{array}{lll}
		\phi_1''+2(u_1^{2}+u_2^{2}+u_3^{2})\phi_1+4(u_1\phi_1+u_2\phi_2+u_3\phi_3)u_1=-\mu_1\phi_1 \ \ \mbox{in}\, \ \R,\\[3mm]
		\phi_2''+2(u_1^{2}+u_2^{2}+u_3^{2})\phi_2+4(u_1\phi_1+u_2\phi_2+u_3\phi_3)u_2=-\mu_2\phi_2 \ \ \mbox{in}\, \ \R,\\[3mm]
		\phi_3''+2(u_1^{2}+u_2^{2}+u_3^{2})\phi_3+4(u_1\phi_1+u_2\phi_2+u_3\phi_3)u_3=-\mu_3\phi_3 \ \ \mbox{in}\, \ \R,
	\end{array}\right.
\end{equation}
where the parameters $\mu_1\leq\mu_2\leq\mu_3<0$ are as in \eqref{GPs}.
The existing works (cf. \cite{A,BS,LW}) show that it is usually necessary to understand the structure of the solution $\phi=(\phi_1, \phi_2, \phi_3)$ for (\ref{Lines}) in $H^1(\R)^3$, in order to analyze the refined limiting behavior of nonlinear Schr\"{o}dinger systems.
Employing Theorem \ref{thm1.1}, we are actually able to prove the following nondegeneracy of the general solution $u$ for \eqref{GPs}.

\begin{thm}\label{thm-nd}
For $\mu_1\leq\mu_2\leq\mu_3<0$, consider the linearized system \eqref{Lines} at any solution $u=(u_1,u_2,u_3)$   of \eqref{GPs} satisfying $u_i\not\equiv 0$ for all $i=1,2,3$. Then  the linearized system \eqref{Lines} at $u$ is non-degenerate, in the sense that the dimension of solutions for the linearized system \eqref{Lines} in $H^1(\R)^3$ is three.
\end{thm}

\begin{rem} \label{rem1.1}
Applying \cite[Section 5.2]{Lewin}, the argument of proving Theorem \ref{thm-nd} yields that any solution $u$ of the system  \eqref{GPN} with $N= 2$ must be non-degenerate. We leave the detailed proof for this case to the interested reader.
\end{rem}

\textcolor{black}{
\begin{rem} \label{rem1.2}
We note from Theorem \ref{thm1.1}  that
\[
u_1(x)=u_1(a_1,a_2,a_3,x),\ u_2(x)=u_2(a_1,a_2,a_3,x), \ u_3(x)=u_3(a_1,a_2,a_3,x)
\]
also depend  on $a_1,a_2,a_3$. If $\mu_1<\mu_2<\mu_3$, then it further yields from Theorem \ref{thm-nd} that the kernel of the linearized system \eqref{Lines} is spanned by
\[
\Big\{(\partial_{a_1} u_1,\partial_{a_1} u_2,\partial_{a_1} u_3),(\partial_{a_2} u_1,\partial_{a_2} u_2,\partial_{a_2} u_3),(\partial_{a_3} u_1,\partial_{a_3} u_2,\partial_{a_3} u_3)\Big\}.
\]
For any $x_0\in\R$, since
\[
u_i(a_1,a_2,a_3,x+x_0)=u_i(a_1e^{\eta_1x_0},a_2e^{\eta_2x_0},a_3e^{\eta_3x_0},x),\,\ j=1,2,3,
\]
we obtain that
\[
\begin{aligned}
(u_1',u_2',u_3')=(\partial_{x} u_1,\partial_{x} u_2,\partial_{x} u_3)=\sum_{j=1}^3\eta_j(\partial_{a_j} u_1,\partial_{a_j} u_2,\partial_{a_j} u_3).
\end{aligned}
\]
\end{rem}
}

Theorem \ref{thm-nd} is another application of Theorem \ref{thm1.1}.
As shown in Section 4, the proof of  Theorem \ref{thm-nd} also makes full use of the integral system (\ref{D-1})--(\ref{D-3}), which is essentially the linearized version of  (\ref{I-1})--(\ref{I-3}) associated to the system \eqref{GPs}.
More precisely, the proof of Theorem \ref{thm-nd} is divided into the following four different cases:
\begin{equation*}\begin{split}
\mbox{\em Case 1.}\ \   \mu_1<\mu_2<\mu_3<0; \quad
\mbox{\em Case 2.}\ \ \mu_1=\mu_2<\mu_3<0;\\
\mbox{\em Case 3.}\ \ \mu_1<\mu_2=\mu_3<0;\quad
\mbox{\em Case 4.}\ \ \mu_1=\mu_2=\mu_3<0 ;
\end{split}
\end{equation*}
Case 1 of Theorem \ref{thm-nd} is proved directly by applying Theorem \ref{thm1.1} and (\ref{D-1})--(\ref{D-3}). In spite of this fact, the argument of addressing Case 1 does not work for the rest three cases. To address Cases 2 and 4 of Theorem \ref{thm-nd}, we shall actually utilize  from Proposition \ref{prop2A} that for any solution $u=(u_1,u_2,u_3)$ of \eqref{GPs}, $\mu_1=\mu_2$ is the first eigenvalue of
 \begin{equation}\label{1:op}
	L_u:= -\frac{d^2}{dx^2}-2(u_1^2+u_2^2+u_3^2)\ \ \mbox{in}\,\ L^2(\R),
\end{equation}
where we also make full use of the fact that the  non-degeneracy of the system  \eqref{GPN} holds essentially for both $N=2$ and $N=1$.

Unfortunately, it is  more challenging to handle with the rest Case 3  of Theorem \ref{thm-nd},  for which case the parameter $\mu_2<0$ is not the first eigenvalue of $L_u$ in $L^2(\R)$. Since $\mu_2=\mu_3<0$ holds for  Case 3, the  eigenvalue $\mu_2$ of  $L_u$ in $L^2(\R)$ may be multiple, which would lead to the degeneracy of solutions for the system \eqref{Lines}. To exclude this possibility, in Section 4 we shall employ the comparison principle to analyze the limiting behavior   of solutions $v\in H^1(\R)$ satisfying $(L_u-\mu_2)v(x)=0$ as $x\to-\infty$, from which we finally conclude that $\mu_2<0$ must be a simple eigenvalue of $L_u$ in $L^2(\R)$. The rest proof of Case 3 is then  similar to that of Case 2.   As a byproduct, the proof of Theorem \ref{thm-nd} implies that for any $\mu_1\leq \mu_2\leq \mu_3\leq 0$, each $\mu_i$ ($i=1,2,3$) of  \eqref{GPs} is always a simple eigenvalue of $L_{u}$ in $L^2(\R)$.

The last section of this paper is devoted to the  $N-$component Schr\"{o}dinger system (\ref{GPN}), where $\mu_1\leq\mu_2\leq\cdots\leq\mu_N<0$ and $N>3$. We shall mainly discuss whether all main results of the present paper can be generalized to the  $N-$component Schr\"{o}dinger system (\ref{GPN}) with any  $N>3$. Roughly speaking, it is proved in Section 5 that the system (\ref{GPN}) still admits the integral system  for any  $N>3$, see Lemma \ref{lem5-1} for more details. This stimulates us to guess in Conjecture 5.1 that the solutions of the system (\ref{GPN})  can be still classified for any  $N>3$, see (\ref{5:N.3}), whose proof is however involved with the more complicated algebraic analysis. Following Conjecture 5.1, one may further guess from Conjecture 5.2 that \eqref{GPN} does not admit any orthonormal solution $u=(u_1,u_2,\cdots,u_N)\in H^1(\R)^N$, where $(u_i, u_j)=\delta_{ij}$ holds for all $1\leq i,j\leq N$ and $N>3$. Of course, motivated by Theorem \ref{thm1.2}, it is our another guess that for any $N\ge 3$, one cannot expect generally the uniqueness  of normalized solutions for the system \eqref{GPN}, which depends on the exact value of $\mu_1\leq\mu_2\leq\cdots\leq\mu_N<0$. Finally, suppose $u=(u_1,u_2,\cdots,u_N)$ is a solution of \eqref{GPN}, then we  consider the following linearized system of \eqref{GPN} around the solution $u$:
\begin{equation}\label{1:GPNLN}
	\begin{array}{lll}
		u_i''+2\Big(\sum_{k=1}^Nu_k^2\Big)\phi_i+4\Big(\sum_{k=1}^Nu_k\phi_k\Big)u_i=-\mu_i\phi_i \ \,\ \mbox{in}\, \ \R , \ \  i=1, 2, \cdots, N.
	\end{array}
\end{equation}
We finally discuss in Section 5 that if Conjecture 5.1 holds true, then we  have the non-degeneracy of (\ref{1:GPNLN}), see those around Theorem \ref{thm5.2} for more details.

Lastly, we mention that the system (\ref{GPN}) belongs to a large class of $N$-component Schr\"{o}dinger systems
\begin{equation}\label{Ncomponent}
	\begin{array}{lll}
		\Delta u_i + \sum_{j=1}^N \beta_{ij} u_j^2 u_i=-\mu_iu_i \ \,\ \mbox{in}\, \ \R^d , \ \  i=1, 2, \cdots, N,
	\end{array}
\end{equation}
where $\beta_{ij}, \mu_i\in \R$, $d\ge 1$ and $N\ge 1$. For the studies of (\ref{Ncomponent}) and its normalized solutions, we refer to \cite{BS, BZZ, CS, LW, WZZ, WW} and the references therein. However, this paper seems to be the first kind for a complete classification on the solutions of (\ref{Ncomponent}) with $N\geq 3$.

This paper is organized as follows. In Section 2, we employ Hirota's bilinearisation method \cite{H} to analyze the general solutions of  \eqref{GPs}, after which we shall classify in Subsection 2.1 the solutions of \eqref{GPs} for the degenerate  cases. In Section 3, we first carry with the further analysis of general solutions for \eqref{GPs}, based on which we shall complete in Subsection 3.1 the proofs of  Theorems \ref{thm1.1} $\&$ \ref{thm1.2} and Corollary \ref{thm1.3}. Section 4 is devoted to the proof of Theorem \ref{thm-nd} on the nondegeneracy of solutions for \eqref{GPs}. Finally, we discuss in Section 5 whether all main results of the present paper can be extended generally to the  $N-$component Schr\"{o}dinger system (\ref{GPN}).

\section{Analysis of General Solutions}

In this section, we first analyze the general solutions $u=(u_1,u_2,u_3)$ of \eqref{GPs}, where $\mu_1\leq\mu_2\leq\mu_3<0$. In Subsection 2.1, we then further discuss the classification on the   solutions of \eqref{GPs} for the degenerate cases where   $\mu_2= \mu_i$ holds for either $i=1$ or $i=3$ or $i=1, 3$.

Stimulated by \cite{Lewin,RL,RSL}, which employ Hirota's bilinearisation method \cite{H}, we set
\begin{equation}\label{1.1}
u_i=\frac{g_i}{f},\quad i=1,2,3.
\end{equation}
It then follows from \eqref{GPs} that $f$ and $g_i$ satisfy
\begin{equation}\label{1-1}
	\left\{\begin{array}{lll}\displaystyle
		fg_i''+g_if''-2f'g_i'+\mu_ifg_i=0,\ \ i=1,2,3,\\[3mm]
		(f')^2-ff''+(g_1^2+g_2^2+g_3^2)=0.
	\end{array}\right.
\end{equation}
Using Hirota's notation, this is of the form:
\[
D(f,g_i)+\mu_i fg_i=0,\ \ D(f, f)=\frac{1}{2}\big(g_1^2+g_2^2+g_3^2\big),\ \ i=1,2,3,
\]
where the bilinear form $D(f,g)$ satisfies $D(f,g)=fg''+f''g-2f'g'$. We consider the following formal expansions:
\begin{equation}\label{1.2}
	g_i=g_{i1}X+g_{i3}X^3+g_{i5}X^5,\quad i=1,2,3,
\end{equation}
and
\begin{equation}\label{1.3}
f=1+f_2X^2+f_4X^4+f_6X^6,
\end{equation}
so that
\begin{equation}\label{1.4}
u_i=\frac{g_{i1}X+g_{i3}X^3+g_{i5}X^5}{1+f_2X^2+f_4X^4+f_6X^6},\quad i=1,2,3.
\end{equation}
We then analyze the following cascade of equations in powers of $X$: for $\eta_i=\sqrt{|\mu_i|}>0,$
\begin{align*}
&g_{i1}''+\mu_ig_{i1}\tag{A-1} = 0,\quad i=1,2,3,\\[2mm]
&-f_2''+(g_{11}^{2}+g_{21}^{2}+g_{31}^{2})\tag{A-2} = 0 ,\\[2mm]
&g_{i3}''+(f_2g_{i1}''+g_{i1}f_2''-2g_{i1}'f_2')+\mu_i(g_{i3}+g_{i1}f_2)\tag{A-3}=0,\quad i=1,2,3,\\[2mm]
&-f_4''+(f_2')^2-f_2f_2''+2(g_{11}g_{13}+g_{21}g_{23}+g_{31}g_{33})\tag{A-4} = 0,\\[2mm]
&g_{i5}''+(f_2g_{i3}''+g_{i3}f_2''-2g_{i3}'f_2')+(f_4g_{i1}''+g_{i1}f_4''-2f_4'g_{i1}')\\
&+\mu_i(g_{i5}+f_2g_{i3}+g_{i3}f_2+f_4g_{i1})\tag{A-5} = 0,\quad i=1,2,3,\\[2mm]
&-f_6''-(f_2f_4''+f_4f_2''-2f_2'f_4')+ \sum_{i=1}^3(2g_{i1}g_{i5}+g_{i3}^2)\tag{A-6}=0,\\[2mm]
&(f_6g_{i1}''+g_{i1}f_6''-2f_6'g_{i1}')+(f_4g_{i3}''+g_{i3}f_4''-2f_4'g_{i3}')\\[2mm]
&+(f_2g_{i5}''+g_{i5}f_2''-2f_2'g_{i5}')+\mu_i(f_6g_{i1}+f_4g_{i3}+f_2g_{i5})\tag{A-7}=0 ,\quad i=1,2,3,\\[2mm]
&(f_4')^2-f_4''f_4+(f_6f_2''+f_2f_6''-2f_2'f_6')+2(g_{13}g_{15}+g_{23}g_{25}+g_{33}g_{35})\tag{A-8}=0,\\[2mm]
&(f_4g_{i5}''+g_{i5}f_4''-2g_{i5}'f_4')+(f_6g_{i3}''+g_{i3}f_6''-2g_{i3}'f_6')\\[2mm]
&+\mu_i(f_4g_{i5}+f_6g_{i3})\tag{A-9}=0,\quad i=1,2,3,\\[2mm]
&-(f_6f_4''+f_4f_6''-2f_4{'}f_6')+ (g^2_{15}+g^2_{25}+g^2_{35})\tag{A-10}=0,\\[2mm]
&(f_6g_{i5}''+g_{i5}f_6''-2g_{i5}'f_6')+\mu_ig_{i5}f_6\tag{A-11}=0,\quad i=1,2,3,\\[2mm]
&(f_6')^2-f_6f_6''\tag{A-12}=0.
\end{align*}
The following lemma gives the formal expansions of (\ref{1.2}) and (\ref{1.3}).

\begin{lem}\label{lem-1}
There exist solutions of (A-1)--(A-12) in the following form:
\begin{equation}\label{s-1}
g_{i1}=a_ie^{\eta_ix}, \quad \eta_i=\sqrt{|\mu_i|}>0,\quad i=1,2,3,
\end{equation}
\begin{equation}\label{s-2}
f_2=\sum_{j=1}^{3}\frac{a_j^2}{4\eta^2_j}e^{2\eta_jx},
\end{equation}
\begin{equation}\label{s-3}
	g_{i3}=\sum_{j=1}^{3}\frac{a_ia_j^2(\eta_i-\eta_j)}{4\eta_j^2(\eta_i+\eta_j)}e^{(2\eta_j+\eta_i)x},\quad  i=1,2,3,
\end{equation}
\begin{equation}\label{s-4}
f_4=\sum_{1\leq j<k\leq 3}\frac{a_j^2a_k^2(\eta_j-\eta_k)^2}{16\eta_j^2\eta_k^2(\eta_j+\eta_k)^2}e^{2(\eta_j+\eta_k)x},
\end{equation}
\begin{equation}\label{s-5}
g_{i5}=\sum_{1\leq j<k\leq 3}\frac{a_ia_j^2a_k^2(\eta_j-\eta_k)^2(\eta_i-\eta_j)(\eta_i-\eta_k)}{16\eta_j^2\eta_k^2(\eta_j+\eta_k)^2(\eta_i+\eta_j)(\eta_i+\eta_k)}e^{(\eta_i+2\eta_j+2\eta_k)x},\quad  i=1,2,3,
\end{equation}
and
\begin{equation}\label{s-6}
f_{6}=\frac{a_1^2a_2^2a_3^2(\eta_1-\eta_2)^2
(\eta_1-\eta_3)^2(\eta_2-\eta_3)^2}{64\eta_1^2\eta_2^2\eta_3^2(\eta_1+\eta_2)^2(\eta_1+\eta_3)^2
(\eta_2+\eta_3)^2}e^{2(\eta_1+\eta_2+\eta_3)x},
\end{equation}
where $a_1, a_2$ and $a_3$ are arbitrary.
\end{lem}

\noindent{\bf Proof.}  We obtain from (A-1) that for $i=1, 2, 3$, $g_{i1}$ satisfies
\begin{equation}\label{26-4}
g_{i1}(x)
= c_{i1} e^{\sqrt{|\mu_i|}x}
+ c_{i2} e^{-\sqrt{|\mu_i|}x},
\quad c_{i1},\, c_{i2}\in\R.
\end{equation}
In the following we construct the solutions of (A-1)--(A-12) by setting $c_{i2}=0$ for (\ref{26-4}). Following this, we specifically choose that $g_{i1}$ satisfies  (\ref{s-1}), where $a_1, a_2$ and $ a_3$  are arbitrary.

Inserting the explicit form \eqref{s-1} into (A-2), we then obtain \eqref{s-2}.
Following the explicit forms \eqref{s-1} and \eqref{s-2}, we have
\[
g_{i3}''+\mu_ig_{i3}+\sum_{j=1}^3\frac{a_ia_j^2}{\eta_j}e^{(\eta_i+2\eta_j)x}=0,
\]
which further yields \eqref{s-3}.
Substituting the explicit forms \eqref{s-1}--\eqref{s-3} into (A-4), it follows that
\[
f_4''-\sum_{1\leq j< k\leq 3}\frac{a_j^2a_k^2(\eta_j-\eta_k)^2}{4\eta_j^2\eta_k^2}e^{2(\eta_j+\eta_k)x}=0,
\]
which implies that
\[
f_4=\sum_{1\leq j< k\leq 3}\frac{a_j^2a_k^2(\eta_j-\eta_k)^2}{16\eta_j^2\eta_k^2(\eta_j+\eta_k)^2}e^{2(\eta_j+\eta_k)x},
\]
$i.e.,$ (\ref{s-4}) holds true.
Similarly, we next substitute \eqref{s-1}--\eqref{s-4} into (A-5), which gives that
\[
g_{i5}''+\mu_ig_{i5}+\sum_{1\leq j< k\leq 3}\frac{a_ia_j^2a_k^2(\eta_i+\eta_j+\eta_k)(\eta_j-\eta_k)^2(\eta_i-\eta_j)
(\eta_k-\eta_i)}{4\eta_j^2\eta_k^2(\eta_j+\eta_k)^2(\eta_i+\eta_j)(\eta_i+\eta_k)}e^{(\eta_i+2\eta_j+2\eta_k)x}=0.
\]
This yields that
\[
g_{i5}=\sum_{1\leq j< k\leq 3}\frac{a_ia_j^2a_k^2(\eta_j-\eta_k)^2(\eta_i-\eta_j)(\eta_i-\eta_k)}{16\eta_j^2
\eta_k^2(\eta_j+\eta_k)^2(\eta_i+\eta_j)(\eta_i+\eta_k)}e^{(\eta_i+2\eta_j+2\eta_k)x},
\]
$i.e.,$ (\ref{s-5}) holds true.

Applying \eqref{s-1}--\eqref{s-5}, we finally deduce from (A-6) that
\begin{equation}\label{2:OKK}
\begin{split}
f_6''&=\big(2f_2'f_4'-f_2f_4''-f_4f_2''\big)+\sum_{i=1}^3(2g_{i1}g_{i5}+g_{i3}^2)\\
&=\sum_{i=1}^3\sum_{j=1}^3\sum_{k=1}^3\frac{a_i^2a_j^2a_k^2}{64\eta_i^2\eta_j^2\eta_k^2}\Big(-\frac{(\eta_j+\eta_k-\eta_i)^2(\eta_j-\eta_k)^2}{(\eta_j+\eta_k)^2}\\
&\quad+\frac{4(\eta_i-\eta_j)(\eta_i-\eta_k)\eta_i^2(\eta_j^2+\eta_k^2)}{(\eta_i+\eta_j)(\eta_i+\eta_k)(\eta_j+\eta_k)^2}\Big)e^{2(\eta_i+\eta_j+\eta_k)x}\\
&:=\sum_{i=1}^3\sum_{j=1}^3\sum_{k=1}^3\frac{a_i^2a_j^2a_k^2}{64\eta_i^2\eta_j^2\eta_k^2}\Big(A_{i,j,k}+B_{i,j,k}\Big)e^{2(\eta_i+\eta_j+\eta_k)x},
\end{split}
\end{equation}
where
\[
A_{i,j,k}:=-\frac{(\eta_j+\eta_k-\eta_i)^2(\eta_j-\eta_k)^2}{(\eta_j+\eta_k)^2},
\]
and
\[
B_{i,j,k}:=\frac{4(\eta_i-\eta_j)(\eta_i-\eta_k)\eta_i^2(\eta_j^2+\eta_k^2)}{(\eta_i+\eta_j)(\eta_i+\eta_k)(\eta_j+\eta_k)^2}.
\]
Note that the term $$\frac{a_i^2a_j^2a_k^2}{64\eta_i^2\eta_j^2\eta_k^2}e^{2(\eta_i+\eta_j+\eta_k)x}$$
is invariant under the permutations of $i,j,k$.

One can check that if $\{i,j,k\}\ne\{1,2,3\}$, then
$$\sum\limits_{\hbox{all permutations of}\,\,i,j,k}(A_{i,j,k}+B_{i,j,k})=0.$$
Therefore, we derive from (\ref{2:OKK}) that
\[
\begin{aligned}
f_6''
&=\frac{a_1^2a_2^2a_3^2}{64\eta_1^2\eta_2^2\eta_3^2}\sum\limits_{ \stackrel{\hbox{all permutations of $i,j,k$,} }{ \hbox{where}\,  \{i, j, k\}=\{1, 2, 3\}}}(A_{i,j,k}+B_{i,j,k})e^{2(\eta_1+\eta_2+\eta_3)x}\\
&=\frac{a_1^2a_2^2a_3^2(\eta_1+\eta_2+\eta_3)^2(\eta_1-\eta_2)^2(\eta_1-\eta_3)^2(\eta_2-\eta_3)^2}{64\eta_1^2\eta_2^2\eta_3^2(\eta_1+\eta_2)^2(\eta_1+\eta_3)^2(\eta_2+\eta_3)^2}e^{2(\eta_1+\eta_2+\eta_3)x},
\end{aligned}
\]
which further implies that
\[
f_6
=\frac{a_1^2a_2^2a_3^2(\eta_1-\eta_2)^2(\eta_1-\eta_3)^2(\eta_2-\eta_3)^2}
{64\eta_1^2\eta_2^2\eta_3^2(\eta_1+\eta_2)^2(\eta_1+\eta_3)^2(\eta_2+\eta_3)^2}e^{2(\eta_1+\eta_2+\eta_3)x}.
\]
This proves (\ref{s-6}), and we are done.
\qed


Applying the symmetry of the system  \eqref{GPs}, we next derive the following three constants of motion for \eqref{GPs}.

\begin{lem}\label{lem2.2} For $\mu_1\leq\mu_2\leq\mu_3<0$, suppose $u=(u_1,u_2,u_3)$ is a  solution of \eqref{GPs} in $H^1(\R)^3$. Then we have the following three constants of motion:	
\begin{equation}\label{I-1}
\begin{split}
(u_1')^2+(u_2')^2+(u_3')^2+(u_1^2+u_2^2+u_3^2)^2+\mu_1u_1^2+\mu_2u_2^2+\mu_3u_3^2=0,	
\end{split}
\end{equation}
\begin{equation}\label{I-2}
	\begin{split}
&(u_1'u_2-u_1u_2')^2+(u_1'u_3-u_1u_3')^2+(u_3'u_2-u_3u_2')^2\\[2mm]
&\,\,+(u_1^2+u_2^2+u_3^2)\big[(\mu_2+\mu_3)u_1^2+(\mu_1+\mu_3)u_2^2+(\mu_1+\mu_2)u_3^2\big]\\[2mm]
&\,\,+\big[(\mu_2+\mu_3)(u'_1)^2+(\mu_1+\mu_3)(u'_2)^2+(\mu_1+\mu_2)(u'_3)^2\big]\\[2mm]
&\,\,+\big[\mu_1(\mu_2+\mu_3)u_1^2+\mu_2(\mu_1+\mu_3)u_2^2+\mu_3(\mu_1+\mu_2)u_3^2\big]=0,
	\end{split}
\end{equation}
and
\begin{equation}\label{I-3}
\begin{split}
& \mu_3(u_1'u_2-u_1u_2')^2+\mu_2(u_1'u_3-u_1u_3')^2+\mu_1(u_3'u_2-u_3u_2')^2\\[2mm]
&\,\,+(u_1^2+u_2^2+u_3^2)( \mu_2\mu_3u_1^2+\mu_1\mu_3u_2^2+\mu_1\mu_2u_3^2)\\[2mm]
&\,\,+\big[\mu_2\mu_3(u'_1)^2+\mu_1\mu_3(u'_2)^2+\mu_1\mu_2(u'_3)^2\big]\\[2mm]
&\,\,+\mu_1\mu_2\mu_3(u_1^2+u_2^2+u_3^2)=0.
\end{split}
\end{equation}
\end{lem}
\noindent{\bf Proof.} Since $(u_1,u_2,u_3)\in H^1(\R)^3$ and $\mu_1\leq\mu_2\leq\mu_3<0$, it is standard to deduce from \eqref{GPs} that both $u_i$ and $u_i'$ vanish  at infinity, where $i=1,2,3$.  Therefore, the above three identities can be derived in a similar way, i.e., by multiplying the
$j$-th equation in \eqref{GPs} by an integrating factor, where $j=1,2,3$, summing up the resulting identities, and finally taking the integration on $\R$.

To derive \eqref{I-1}, more precisely, we multiply the $j$-th equation of \eqref{GPs} by $u_j'$, where $j=1,2,3$.  Then we add together the resulting identities and integrate it on $\R$, where we use the fact that both $u_j$ and $u_j'$ vanish at infinity. This thus gives \eqref{I-1}.

To obtain \eqref{I-2}, we multiply the $j$-th equation of \eqref{GPs} by
$$\sum_{1\leq k\leq 3,\,k\neq j }\big[2u_k(u_j'u_k-u_ju_k')+2\mu_ku_j'\big],\ \ j=1,2,3.$$
  We then add together the resulting identities and integrate it on $\R$ again, which thus yields \eqref{I-2}.

Multiply  the $j$-th equation of \eqref{GPs} by
$$\sum_{1\leq i<k\leq 3,\, i\neq j,\, k\neq j }\big[2\mu_ku_i(u_j'u_i-u_ju_i')+2\mu_iu_k(u_j'u_k-u_ju_k')+2\mu_i\mu_ku_j'\big], \ \ j=1,2,3,$$
which then gives  \eqref{I-3} in the similar way as above.  This therefore completes  the proof of the lemma.   \qed

\subsection{Classification of solutions for degenerate  cases}
In this subsection, we classify the solutions of \eqref{GPs} for the degenerate  cases where   $\mu_2= \mu_i$ holds for either $i=1$ or $i=3$ or $i=1, 3$. Our main results of this subsection can be stated as the following proposition.

\begin{prop}\label{prop2A}
For $\mu_1\leq \mu_2\leq \mu_3<0$,  let $u=(u_1,u_2,u_3)\in H^1(\R)^3$ be a  solution of \eqref{GPs}, and denote $\eta_i=\sqrt{|\mu_i|}>0$ for $ i=1, 2, 3$. Then we have the following conclusions:
\begin{itemize}
\item[(i.)] If $\mu_1=\mu_2<\mu_3<0$,  then there exists $(a_1,a_2,a_3)\in\R^3$ such that $u=(u_1,u_2,u_3)$ satisfies
\begin{equation}\label{c1}
	\left\{\begin{array}{lll}
u_1(x)=\frac{a_1e^{\eta_1x}}{f(x)}\Big[1+\frac{a_3^2(\eta_1-\eta_3)}{4\eta_3^2(\eta_1+\eta_3)}e^{2\eta_3x}\Big],\\[3mm]
u_2(x)=\frac{a_2e^{\eta_1x}}{f(x)}\Big[1+\frac{a_3^2(\eta_1-\eta_3)}{4\eta_3^2(\eta_1+\eta_3)}e^{2\eta_3x}\Big],\\[3mm]
u_3(x)=\frac{a_3e^{\eta_3x}}{f(x)}\Big[1+\frac{(a_1^2+a_2^2)(\eta_3-\eta_1)}{4\eta_1^2(\eta_1+\eta_3)}e^{2\eta_1x}\Big],
	\end{array}\right.
\end{equation}
where $f$ is given by
\[
f(x)=1+\frac{a_3^2}{4\eta_3^2}e^{2\eta_3x}+\frac{(a_1^2+a_2^2)}{4\eta_1^2}e^{2\eta_1x}+\frac{(a_1^2+a_2^2)a_3^2(\eta_3-\eta_1)^2}{16\eta_1^2\eta_3^2(\eta_1+\eta_3)^2}e^{2(\eta_1+\eta_3)x}.
\]
Furthermore, if $u^2_1+u^2_2\not\equiv 0$ and $u_3\not\equiv 0$, then we have $\int_{\R}(u^2_1+u^2_2)dx=2\eta_1$ and $\int_{\R}u^2_3dx=2\eta_3$.

\item[(ii.)]  If $\mu_1<\mu_2=\mu_3<0$,  then there exists $(a_1,a_2,a_3)\in\R^3$ such that $u=(u_1,u_2,u_3)$ satisfies
\begin{equation}\label{c2}
	\left\{\begin{array}{lll}
		 u_1(x)=\frac{a_1e^{\eta_1x}}{f(x)}\Big[1+\frac{(a_2^2+a_3^2)(\eta_1-\eta_3)}{4\eta_3^2(\eta_1+\eta_3)}e^{2\eta_3x}\Big],\\[3mm]
		 u_2(x)=\frac{a_2e^{\eta_2x}}{f(x)}\Big[1+\frac{a_1^2(\eta_2-\eta_1)}{4\eta_1^2(\eta_1+\eta_2)}e^{2\eta_1x}\Big],\\[3mm]
		u_3(x)=\frac{a_3e^{\eta_2x}}{f(x)}\Big[1+\frac{a_1^2(\eta_2-\eta_1)}{4\eta_1^2(\eta_1+\eta_2)}e^{2\eta_1x}\Big],
	\end{array}\right.
\end{equation}
where $f$ is given by
\[
f(x)=1+\frac{a_1^2}{4\eta_1^2}e^{2\eta_1x}+\frac{(a_2^2+a_3^2)}{4\eta_3^2}e^{2\eta_3x}+\frac{a_1^2(a_2^2+a_3^2)(\eta_1-\eta_3)^2}{16\eta_1^2\eta_3^2(\eta_1+\eta_3)^2}e^{2(\eta_1+\eta_3)x}.
\]
Furthermore, if $u_1\not\equiv 0$ and $u^2_2+u^2_3\not\equiv 0$, then we have $\int_{\R}u^2_1dx=2\eta_1$ and $\int_{\R}(u^2_2+u^2_3)dx=2\eta_2$.

\item[(iii.)] If $\mu_1=\mu_2=\mu_3<0$,  then there exists $(a_1,a_2,a_3)\in\R^3$ such that $u=(u_1,u_2,u_3)$ satisfies
\begin{equation}\label{c23}
u_i(x)=\frac{a_ie^{\eta_1x}}{f(x)},\,\ i=1,2,3,
\end{equation}
where $f(x)$ is given by
\[
f(x)=1+\frac{(a_1^2+a_2^2+a_3^2)}{4\eta_1^2}e^{2\eta_1x}.
\]
Furthermore,  if $\sum_{i=1}^3u^2_i\not\equiv 0$, then we have $\sum_{i=1}^3\int_{\R}u^2_idx=2\eta_1$.

\end{itemize}
\end{prop}

\noindent{\bf Proof.} (i). Suppose $\mu_1=\mu_2<\mu_3<0$, then we have $\mu_3>\mu_0$, where $\mu_0<0$ is the lowest eigenvalue of $-\frac{d^2}{dx^2}-2(u_1^2+u_2^2+u_3^2) $ in $\R$. Thus, $u_3$ vanishes at some point, and hence we can suppose $u_3(0)=0$ by a suitable translation. Also, we may assume that  $u_3'(0)\neq0$, since otherwise we have $u_3\equiv 0$ in $\R$, and a similar argument of \cite[Lemma 15]{Lewin} further gives (\ref{c1}). Similarly, we may assume that $u_2\not\equiv 0$, since otherwise the classification (\ref{c1}) follows directly from \cite[Lemma 15]{Lewin}.

Since $\mu_1=\mu_2<\mu_3<0$ and $u_2\not\equiv 0$, we now claim that
\begin{equation}\label{clm1}
u_1\equiv ku_2\ \ \hbox{holds for some $k\in\R$}.
\end{equation}
Indeed, because $u_3(0)=0$ and $\mu_1=\mu_2$, we derive from \eqref{I-1} and \eqref{I-2} that
\begin{equation}\label{Ns-21A}
	\left\{\begin{array}{lll}
u_1'(0)^2+u_2'(0)^2+u_3'(0)^2=-\big[(u_1(0)^2+u_2(0)^2)^2+\mu_1(u_1(0)^2+u_2(0)^2)\big],\\[3mm]
		\quad \big[u_1'(0)u_2(0)-u_1(0)u_2'(0)\big]^2+(\mu_1+\mu_3)(u_1'(0)^2+u_2'(0)^2)\\[3mm]
		\quad+\big[u_1(0)^2+u_2(0)^2+2\mu_1\big]u_3'(0)^2\\[3mm]
		=-(\mu_1+\mu_3)[u_1(0)^2+u_2(0)^2]^2-\mu_1(\mu_1+\mu_3)(u_1(0)^2+u_2(0)^2).
	\end{array}\right.
\end{equation}
\textcolor{black}{Calculating
$\eqref{I-3}-\mu_3\times\eqref{I-2}+\mu_3^2\times\eqref{I-1}$, it  yields that}
\[
\big[(u_1(0)^2+u_2(0)^2)-(\mu_3-\mu_1)\big]u_3'(0)^2=0.
\]
Since $u_3'(0)^2\neq 0$,  the above three equations give that
\begin{equation}\label{Ns-1}
	\left\{\begin{array}{lll}
u_1'(0)^2+u_2'(0)^2+u_3'(0)^2=-\big[\big(u_1(0)^2+u_2(0)^2\big)^2+\mu_1u_1(0)^2+\mu_2u_2(0)^2\big],\\[3mm]
		\big[u_1'(0)u_2(0)-u_1(0)u_2'(0)\big]^2
		+\big[u_1(0)^2+u_2(0)^2+(\mu_1-\mu_3)\big]u_3'(0)^2=0,\\[3mm]
	u_1(0)^2+u_2(0)^2=(\mu_3-\mu_1),
	\end{array}\right.
\end{equation}
where the second equation is derived   \textcolor{black}{by  the calculations of $(\mu_1+\mu_3)\times\eqref{Ns-21A}_1-\eqref{Ns-21A}_2$.}

By the last two equations of (\ref{Ns-1}), it yields that
\begin{equation}\label{Sep3}
 u_1'(0)u_2(0)-u_1(0)u_2'(0)=0.
\end{equation}
The above equation implies that the  initial values $(u_1(0),u_1'(0))$ and $(u_2(0),u_2'(0))$ of $u_1$, $u_2$ are proportional. Combined with the fact that $u_1$, $u_2$ satisfy the same second order ODE:
\[
u''+2(u_1^2+u_2^2+u_3^2)u+\mu_1u=0,
\]
we get from the uniqueness argument of ODE that the claim \eqref{clm1} holds true.

Applying \eqref{clm1} and setting
\[
\tilde u=\sqrt{1+k^2}u_2,
\] then we derive from \eqref{GPs} that $(\tilde{u},u_3)$ satisfies the following two component system:
\begin{equation}\label{eq2}
\left\{\begin{array}{lll}
\tilde{u}''+2\big[\tilde{u}^2+u_3^2\big]\tilde{u}+\mu_1\tilde{u}=0,\\[2mm]
u_3''+2\big[\tilde{u}^2+u_3^2\big]u_3+\mu_3u_3=0.
\end{array}\right.
\end{equation}
Then \cite[Lemma 15]{Lewin} implies that there exists $(\tilde{a},a_3)$ such that
$(\tilde{u},u_3)$ satisfies
\begin{equation}\label{O1}
\left\{\begin{array}{lll}
\tilde{u}(x)=\frac{\tilde{a}e^{\eta_1x}}{f(x)}\Big[1+\frac{a_3^2(\eta_1-\eta_3)}{4\eta_3^2(\eta_1+\eta_3)}e^{2\eta_3x}\Big],\\[3mm]
u_3(x)=\frac{a_3e^{\eta_3x}}{f(x)}\Big[1+\frac{\tilde{a}^2(\eta_3-\eta_1)}{4\eta_1^2(\eta_1+\eta_3)}e^{2\eta_1x}\Big],
\end{array}\right.
\end{equation}
where $f$ is given by
\[
f(x)=1+\frac{a_3^2}{4\eta_3^2}e^{2\eta_3x}+\frac{\tilde{a}^2}{4\eta_1^2}e^{2\eta_1x}+\frac{\tilde{a}^2a_3^2(\eta_3-\eta_1)^2}{16\eta_1^2\eta_3^2(\eta_1+\eta_3)^2}e^{2(\eta_1+\eta_3)x}.
\]
Following \eqref{clm1} and \eqref{O1}, we obtain that \eqref{c1} holds true by setting
\[
a_1=\frac{|k\tilde{a}|}{1+k^2}\hbox{sign}\,u_1(0)\,\ \hbox{and}\,\ a_2=\frac{|\tilde{a}|}{1+k^2}\hbox{sign}\,u_2(0).
\]

Finally, if $u_1^2+u_2^2\not\equiv 0$ and $u_3^2\not\equiv 0$, we obtain from \cite[Lemma 16]{Lewin} that
\[
\inte(u_1^2+u_2^2)dx=\inte\tilde{u}^2dx=-\inte\left(\frac{\frac{a_3^2\eta_1}{2\eta_3^2}e^{2\eta_3x}+2\eta_1}{f(x)}\right)'=2\eta_1,
\]
and
\[
\int_{\R}u_3^2dx=-\inte\left(\frac{\frac{\tilde{a}^2\eta_3}{2\eta_1^2}e^{2\eta_1x}+2\eta_3}{f(x)}\right)'=2\eta_3.
\]
This completes the proof of (i).

(ii). Suppose $\mu_1<\mu_2=\mu_3<0$, then we have $\mu_3>\mu_1\geq \mu_0$, where $\mu_0<0$ is the lowest eigenvalue of $-\frac{d^2}{dx^2}-2(u_1^2+u_2^2+u_3^2) $ in $\R$. Thus, $u_3$ vanishes at some point. As before, we may assume that $u_3(0)=0$ and $u_3'(0)\neq0$, which then implies  that $u_3(x)\not \equiv 0$.

\textcolor{black}{Calculating $\eqref{I-3}-\mu_3\times \eqref{I-2}+\mu_3^2\times\eqref{I-1}$, it yields that}
\begin{equation}\label{Sep2}
u_2(0)^2u_3'(0)^2=0,
\end{equation}
which further implies that $u_2(0)=0$, since $u_3'(0)\neq0$. We hence have
$u_2(0)=u_3(0)=0$. Because $u_3'(0)\neq0$, there exists some $k\in\R$ such that $u_2'(0)=ku_3'(0)$.
Since $\mu_1<\mu_2=\mu_3<0$ and $u_3(x)\not \equiv 0$,  the same argument of (\ref{clm1}) then yields that
\begin{equation}\label{clm2}
	u_2\equiv ku_3\,\ \hbox{holds for the above $k\in\R$}.
\end{equation}

Applying \eqref{clm2} and setting
\[
\hat u=\sqrt{1+k^2}u_3,
\]
then we derive from \eqref{GPs} that $(u_1,\hat{u})$ satisfies the following two component system:
\begin{equation}\label{O2}
	\left\{\begin{array}{lll}
u_1''+2\big[u_1^2+\hat{u}^2\big]u_1+\mu_1u_1=0,\\[2mm]	\hat{u}''+2\big[u_1^2+\hat{u}^2\big]\hat{u}+\mu_2\hat{u}=0.
\end{array}\right.
\end{equation}
Using \cite[Lemmas 15 and 16]{Lewin}, the same argument in (i) applies to \eqref{O2} yields that (ii) holds true.

(iii). Suppose $\mu_1=\mu_2=\mu_3<0$. \textcolor{black}{Calculating
$\eqref{I-3}-2\mu_1\times\eqref{I-2}$, it gives that}
\begin{equation*}
\begin{split}
	&\big(u_1'(0)u_2(0)-u_1(0)u_2'(0)\big)^2+\big(u_1'(0)u_3(0)-u_1(0)u_3'(0)\big)^2\\
&+\big(u_3'(0)u_2(0)-u_3(0)u_2'(0)\big)^2=0,
\end{split}
\end{equation*}
which further implies that
\begin{equation}\label{sp-9}
\begin{split}
	u_1'(0)u_2(0)-u_1(0)u_2'(0)&=u_1'(0)u_3(0)-u_1(0)u_3'(0)\\
&=u_3'(0)u_2(0)-u_3(0)u_2'(0)=0.
\end{split}
\end{equation}
Without loss of generality, we may assume that $u_3(x)\not\equiv 0$. It then follows from \eqref{sp-9} that there exist $k_1,k_2\in\R$ such that
\[
u_1(0)=k_1u_3(0),\,\ u_1'(0)=k_1u_3'(0),
\]
and
\[
u_2(0)=k_2u_3(0),\,\ u_2'(0)=k_2u_3'(0).
\]
The same argument of \eqref{clm1} gives that $u_1\equiv k_1u_3$ and $u_2\equiv k_2u_3$.

Thus, $u_3$ satisfies the following equation
\begin{equation}\label{26-3}
u_3''+2(1+k_1^2+k_2^2)u_3^3=-\mu_1u_3,	
\end{equation}
which admits the following constant of motion:
\[
(u_3')^2+(1+k_1^2+k_2^2)u_3^4+\mu_1u_3^2=0.
\]
Since $u_3\not\equiv 0$, the uniqueness of  solutions for \eqref{26-3} implies that $u_3$ is either positive or negative, which further yields that $\mu_1=\mu_2=\mu_3$ is the first eigenvalue of the operator $$-\frac{d^2}{dx^2}-2(u_1^2+u_2^2+u_3^2)\,\ \hbox{in}\  L^2(\R).$$
Hence, \eqref{c23} holds true, where $\frac{1}{1+k_1^2+k_2^2}=\frac{a_3}{a_1^2+a_2^2+a_3^2}$, and $\frac{k_i}{1+k_1^2+k_2^2}=\frac{a_i}{a_1^2+a_2^2+a_3^2}$ holds for $i=1,2$. Direct computations further give that $\sum_{i=1}^3\int_{\R}u_i^2dx=2\eta_1$.
This completes the proof of Proposition \ref{prop2A}.  \qed

The above proof of Proposition \ref{prop2A} shows that if $\mu_2= \mu_i$ holds for either $i=1$ or $i=3$ or $i=1, 3$, then the classification of general solutions   for the system \eqref{GPs} follows essentially from the argument of \cite[Lemma 15]{Lewin}, which unfortunately  is not applied to the general case $\mu_1<\mu_2<\mu_3<0$.

\section{Classification of General Solutions}
In this section we further analyze the properties of general solutions for \eqref{GPs} in the case where $\mu_1<\mu_2<\mu_3<0$, based on which we shall complete in Subsection 3.1 the proofs of Theorems \ref{thm1.1} $\&$ \ref{thm1.2} and  Corollary \ref{thm1.3}.  \textcolor{black}{Up to a suitable translation, we may assume that $u_3(0)=0$ in the whole section. Indeed, if either $u_3>0$ or $u_3<0$, then $\mu_3$ is the first eigenvalue of the operator $-\frac{d^2}{dx^2}-2(u_1^2+u_2^2+u_3^2)\ $ in $L^2(\R)$, which thus implies that $\mu_3\leq\mu_2\leq\mu_1$, a contradiction. }

We start with the following lemma.

\begin{lem}\label{lem:3.1} For $\mu_1<\mu_2<\mu_3<0$, suppose $u=(u_1,u_2,u_3)\in H^1(\R)^3$ is a  solution of \eqref{GPs} satisfying $u_3(x)\not\equiv 0$ in $\R$. Then we have the following conclusions:
	\begin{itemize}
		\item[(i).]  If $u_i(0)=u_3(0)=0$ holds for either $i=1$ or $2$, then we have $u_i\equiv 0$.
		
		\item[(ii).]  If  $u'_1(0)=u'_2(0)=u_3(0)=0$, then we have either $u_1\equiv 0$ or $u_2\equiv 0$.
	\end{itemize}
\end{lem}

\noindent{\bf Proof.}  By the assumption $u_3(0)=0$ of Lemma \ref{lem:3.1}, \textcolor{black}{we derive from
	$\eqref{I-3}-\mu_3\times\eqref{I-2}+\mu_3^2\times\eqref{I-1}$  that}
\begin{equation}\label{Sep7}
(\mu_3-\mu_2)u_1(0)^2+(\mu_3-\mu_1)u_2(0)^2=(\mu_3-\mu_1)(\mu_3-\mu_2).
\end{equation}
We next follow (\ref{Sep7}) to complete the proof.

(i). We only consider the case $i=1$, so that $u_1(0)=u_3(0)=0$, since the other case can be proved similarly. For this case, we first obtain from \eqref{Sep7} that $u_2(0)^2=\mu_3-\mu_2$. Moreover, we derive from \eqref{I-1} and \eqref{I-2} that
\begin{equation}\label{Sep9}
	\left\{\begin{array}{lll}
u_1'(0)^2+u_2'(0)^2+u_3'(0)^2=-\mu_3(\mu_3-\mu_2),\\[3mm]
u_2'(0)^2+u_3'(0)^2=-\mu_3(\mu_3-\mu_2),
	\end{array}\right.
\end{equation}
which implies that $u_1'(0)=0$. Hence $u_1\equiv0$, and we are done.

(ii). Suppose  $u'_1(0)=u'_2(0)=u_3(0)=0$, then we derive from \eqref{I-1}, \eqref{I-2} and \eqref{Sep7} that
\begin{equation}\label{Sep10}
	\left\{\begin{array}{lll}
u_3'(0)^2=-\big[(u_1(0)^2+u_2(0)^2)^2+\mu_1u_1(0)^2+\mu_2u_2(0)^2\big],\\[3mm]
\big[u_1(0)^2+u_2(0)^2+(\mu_2-\mu_3)\big]u_3'(0)^2\\[3mm]
 =-(\mu_2-\mu_1)\big[u_1(0)^2+u_2(0)^2+\mu_1\big]u_1(0)^2,\\[3mm]
(\mu_3-\mu_2)u_1(0)^2+(\mu_3-\mu_1)u_2(0)^2=(\mu_3-\mu_1)(\mu_3-\mu_2).
\end{array}\right.
\end{equation}
By eliminating the variable $u_3'(0)^2$ from the first two equations of \eqref{Sep10}, we obtain that
\begin{equation}\label{Sep6}
	\left\{\begin{array}{lll}
\big[(u_1(0)^2+u_2(0)^2)^2+\mu_1u_1(0)^2+\mu_2u_2(0)^2\big]\big[u_1(0)^2+u_2(0)^2+(\mu_2-\mu_3)\big]\\[3mm]
=(\mu_2-\mu_1)\big[u_1(0)^2+u_2(0)^2+\mu_1\big]u_1(0)^2,\\[3mm]
(\mu_3-\mu_2)u_1(0)^2+(\mu_3-\mu_1)u_2(0)^2=(\mu_3-\mu_1)(\mu_3-\mu_2).\\[3mm]
	\end{array}\right.
\end{equation}
One can verify that the solution $(u_1(0)^2,u_2(0)^2)$ of the above system    has the following three different forms:
$$(u_1(0)^2,u_2(0)^2)=(0,\mu_3-\mu_2),\ \ (u_1(0)^2,u_2(0)^2)=(\mu_3-\mu_1,0),$$ $$(u_1(0)^2,u_2(0)^2)=\Big(-\frac{(\mu_1-\mu_3)\mu_1}{\mu_1-\mu_2},\frac{(\mu_2-\mu_3)\mu_2}{\mu_1-\mu_2}\Big).$$
Here the second component of the last one is negative, which is impossible. Therefore, we have either $u_1(0)=0$, which implies that $u_1(x)\equiv 0$, or $u_2(0)=0$, which implies that $u_2(x)\equiv 0$. Lemma \ref{lem:3.1} is therefore proved.  \qed

In the following, we consider  the last case where $\mu_1<\mu_2<\mu_3<0$, and
\begin{equation}\label{Sep11}
u_1(0)\neq 0,\,\ u_2(0)\neq 0,\,\ u_i'(0)\neq 0\,\ \hbox{for some $1\leq i\leq 2$,\,\ }u_3(0)=0,\,\ u_3'(0)\neq 0.
\end{equation}
Without loss of generality,   the assumption of (\ref{Sep11}) is equivalent to the following form: there exist a constant $p\in\R$ and a nonzero constant $ q\in \R$ such that
\begin{equation}\label{N-4}
	u_3(0)=0,\,\, \frac{u_2(0)}{u_1(0)}=q\not= 0,\,\,  u'_2(0)=pu'_1(0),\ \  u_1'(0)\neq 0,\,\   u_3'(0)\neq 0.
\end{equation}
Define
\begin{equation}\label{S-8}
\begin{split}
f(p):&=-\Big\{\frac{(1+q^2)^2(\mu_3-\mu_1)(\mu_3-\mu_2)}{(\mu_3-\mu_2)+q^2(\mu_3-\mu_1)}+\mu_1+q^2\mu_2\Big\}\\
&\quad-(1+p^2)\frac{q^2(\mu_1-\mu_2)^2(\mu_1^2q^2-\mu_1\mu_3q^2+\mu_2^2-\mu_2\mu_3)}
{[(\mu_3-\mu_2)+q^2(\mu_3-\mu_1)][(\mu_3-\mu_2)p-(\mu_3-\mu_1)q]^2}.
	\end{split}	
\end{equation}
We shall prove in \eqref{SS-9} below that for any given $q\not= 0$, $f(p)$ admits exactly two different zero points denoted by $\underline{p}=\underline{p}(q)<\overline{p}=\overline{p}(q)$, $i.e.$,
\begin{equation}\label{S-9}
	f(\underline{p})=f(\overline{p})=0,  \ \  \hbox{where}\ \  \underline{p}<\overline{p}.
\end{equation}
We now address the following classification for the last case  (\ref{N-4}):

\begin{prop}\label{prop2B} For $\mu_1<\mu_2<\mu_3<0$, suppose the function $u=(u_1,u_2,u_3)\in H^1(\R)^3$ satisfies \eqref{N-4} at $x=0$.
Then
$$\big(u_1(0),u_2(0),u_3(0),u'_1(0),u'_2(0),u'_3(0)\big)$$ satisfies \eqref{I-1}--\eqref{I-3}  at $x=0$, if and only if
\begin{equation}\label{N-17}
q\not =0, \ \ p\in(\underline{p},\,\,\overline{p}),
\end{equation}
where the constants $\underline{p}=\underline{p}(q)$ and $\overline{p}=\overline{p}(q)$ are as in (\ref{S-9}). Moreover, if $p$ satisfies \eqref{N-17}, then \eqref{I-1}--\eqref{I-3}  at $x=0$  admits at most the following four solutions:
\begin{equation}\label{N-21}
	\begin{split}
		\big(u_1(0),u_2(0),0,u'_1(0),u'_2(0),u'_3(0)\big),\qquad\\
\big(-u_1(0),-u_2(0),0,-u'_1(0),-u'_2(0),u'_3(0)\big),\\
		\big(u_1(0),u_2(0),0,u'_1(0),u'_2(0),-u'_3(0)\big),\quad\ \\
\big(-u_1(0),-u_2(0),0,-u'_1(0),-u'_2(0),-u'_3(0)\big).
	\end{split}
\end{equation}
\end{prop}

\noindent{\bf Proof.} Suppose $u=(u_1,u_2,u_3)\in H^1(\R)^3$ is a  solution of \eqref{I-1}--\eqref{I-3} satisfying \eqref{N-4} at $x=0$, then the standard elliptic regularity theory gives that $u\in C^1(\R)^3$.
Applying \eqref{N-4}, we derive from \eqref{I-1} and \eqref{I-2} that
\begin{equation}\label{N-21A}
	\left\{\begin{array}{lll}
(1+p^2)u_1'(0)^2+u_3'(0)^2=-\big[(u_1(0)^2+u_2(0)^2)^2+\mu_1u_1(0)^2+\mu_2u_2(0)^2\big],\\[3mm]
\quad\big[(p-q)^2u_1(0)^2+(\mu_2+\mu_3)+p^2(\mu_1+\mu_3)\big]u_1'(0)^2\\[3mm]
\quad+\big[u_1(0)^2+u_2(0)^2+(\mu_1+\mu_2)\big]u_3'(0)^2\\[3mm]
=-\big[u_1(0)^2+u_2(0)^2\big]\big[(\mu_2+\mu_3)u_1(0)^2+(\mu_1+\mu_3)u_2(0)^2\big]\\[3mm]
\quad -\big[\mu_1(\mu_2+\mu_3)u_1(0)^2+\mu_2(\mu_1+\mu_3)u_2(0)^2\big].
	\end{array}\right.
\end{equation}
\textcolor{black}{Calculating
	$\eqref{I-3}-\mu_3\times\eqref{I-2}+\mu_3^2\times\eqref{I-1}$, it  yields that}
\[
\big[(\mu_3-\mu_2)u_1(0)^2+(\mu_3-\mu_1)u_2(0)^2-(\mu_3-\mu_1)(\mu_3-\mu_2)\big]u_3'(0)^2=0,
\]
where $u_3'(0) \neq 0$ in view of (\ref{N-4}).
Following the above three equations, we then deduce from (\ref{N-4}) that
\begin{equation}\label{N-1}
	\left\{\begin{array}{lll}
		(1+p^2)u_1'(0)^2+u_3'(0)^2=-\big[(1+q^2)^2u_1(0)^2+\mu_1+q^2\mu_2\big]u_1(0)^2,\\[3mm]
		\big[(p-q)^2u_1(0)^2-(\mu_1-\mu_2)\big]u_1'(0)^2
		+\big[(1+q^2)u_1(0)^2+(\mu_2-\mu_3)\big]u_3'(0)^2\\[3mm]
		\quad=(1+q^2)(\mu_1-\mu_2)u_1(0)^4+\mu_1(\mu_1-\mu_2)u_1(0)^2,\\[3mm]
		\big[(\mu_3-\mu_2)+q^2(\mu_3-\mu_1)\big]u_1(0)^2=(\mu_3-\mu_1)(\mu_3-\mu_2),
	\end{array}\right.
\end{equation}
where the constants $q$ and $p$ are as in (\ref{N-4}), and  $\eqref{N-1}_2$ is derived by $\eqref{N-21A}_2-(\mu_1+\mu_3)\times \eqref{N-21A}_1$.
Note  that the values of $u_1(0)^2$ and $u_2(0)^2=q^2u_1(0)^2$ are determined by the last equation of (\ref{N-1}).

By  Cramer's Rule, we obtain from  (\ref{N-1}) that
\[
u_1'(0)^2=\frac{D_1}{D},
\]
where
\[
D=\det\begin{vmatrix}
(1+p^2)& 1\\[3mm]
(p-q)^2u_1(0)^2-(\mu_1-\mu_2)&(1+q^2)u_1(0)^2+(\mu_2-\mu_3)
\end{vmatrix},
\]
and
\[
D_1=\det\begin{vmatrix}
-[(1+q^2)^2u_1(0)^2+\mu_1+q^2\mu_2]u_1(0)^2& 1\\[3mm]
(1+q^2)(\mu_1-\mu_2)u_1(0)^4+\mu_1(\mu_1-\mu_2)u_1(0)^2& (1+q^2)u_1(0)^2+(\mu_2-\mu_3)
\end{vmatrix}.
\]
%
Direct calculations then yield from  (\ref{N-1}) that
\begin{equation*}
\begin{split}
D&=	(pq+1)^2u_1(0)^2-\big[\mu_3-\mu_1+p^2(\mu_3-\mu_2)\big]\\
&=\big[q^2u_1(0)^2-(\mu_3-\mu_2)\big]p^2+2pqu_1(0)^2+\big[u_1(0)^2-(\mu_3-\mu_1)\big]\\
&=-\frac{1}{(\mu_3-\mu_2)+q^2(\mu_3-\mu_1)}\big[(\mu_3-\mu_2)p-(\mu_3-\mu_1)q\big]^2,
\end{split}	
\end{equation*}
and
\begin{equation*}
	\begin{split}
\frac{D_1}{u_1(0)^2}
&=(1+q^2)(\mu_1-\mu_2)u_1(0)^2+\mu_1(\mu_1-\mu_2)+(1+q^2)^3u_1(0)^4\\
&\quad+(\mu_2-\mu_3)(1+q^2)^2u_1(0)^2+(\mu_1+\mu_2q^2)(1+q^2)u_1(0)^2\\
&\quad+(\mu_1+\mu_2q^2)(\mu_2-\mu_3)\\
&=-\frac{q^2 (\mu_1 - \mu_2)^2 (\mu_1^2 q^2 - \mu_1 \mu_3 q^2 + \mu_2^2 - \mu_2 \mu_3)}{ [(\mu_3-\mu_2)+q^2(\mu_3-\mu_1)]^2}.
	\end{split}	
\end{equation*}
We then have
\begin{equation}\label{N-2}
	\begin{split}
0\leq u_1'(0)^2&=\frac{D_1}{D}\\
&=\frac{q^2(\mu_1-\mu_2)^2(\mu_3-\mu_1)(\mu_3-\mu_2)(\mu_1^2q^2-\mu_1\mu_3q^2+\mu_2^2-\mu_2\mu_3)}
{[(\mu_3-\mu_2)+q^2(\mu_3-\mu_1)]^2[(\mu_3-\mu_2)p-(\mu_3-\mu_1)q]^2}.
	\end{split}	
\end{equation}
Since $\frac{D_1}{D}$ is positive for any $p\in\R $ and $0\neq q\in\R$,   there is no restriction on  $p$ and $q$ to ensure that $ u_1'(0)^2> 0$.

We next analyze $u_3'(0)^2$ as follows. It follows from \eqref{N-1} and \eqref{N-2} that for any $q\not =0,$
\begin{align}
0\leq \frac{u_3'(0)^2}{u_1(0)^2}&=-\big[(1+q^2)^2u_1(0)^2+\mu_1+q^2\mu_2\big]\nonumber\\
&\quad-(1+p^2)\frac{q^2(\mu_1-\mu_2)^2(\mu_1^2q^2-\mu_1\mu_3q^2+\mu_2^2-\mu_2\mu_3)}
{[(\mu_3-\mu_2)+q^2(\mu_3-\mu_1)][(\mu_3-\mu_2)p-(\mu_3-\mu_1)q]^2}\nonumber\\
&=-\Big\{\frac{(1+q^2)^2(\mu_3-\mu_1)(\mu_3-\mu_2)}{(\mu_3-\mu_2)+q^2(\mu_3-\mu_1)}+\mu_1+q^2\mu_2\Big\}\label{S-3}\\
&\quad-\frac{q^2(\mu_1-\mu_2)^2(\mu_1^2q^2-\mu_1\mu_3q^2+\mu_2^2-\mu_2\mu_3)}
{[(\mu_3-\mu_2)+q^2(\mu_3-\mu_1)]}\frac{1+p^2}{[(\mu_3-\mu_2)p-(\mu_3-\mu_1)q]^2}\nonumber\\
&=:f(p).\nonumber
\end{align}	
For any $q\not =0,$ since $\mu_1<\mu_2<\mu_3<0$, the second term of $f(p)$ is strictly negative, and
\begin{equation}\label{3C:S-3}
\begin{aligned}
\frac{1+p^2}{[(\mu_3-\mu_2)p-(\mu_3-\mu_1)q]^2}
=&\frac{1}{\big[(\mu_3-\mu_2)\frac{p}{\sqrt{1+p^2}}-(\mu_3-\mu_1)q\frac{1}{\sqrt{1+p^2}}\big]^2}\\
=&\frac{1}{[(\mu_3-\mu_2)^2+(\mu_3-\mu_1)^2q^2]\cos^2\theta_p},
\end{aligned}
\end{equation}
where $\theta_p\in [0, \frac{\pi}{2})\cup (\frac{\pi}{2}, \pi]$ denotes the angle between the vectors $(\mu_3-\mu_2,-(\mu_3-\mu_1)q)$ and $(\frac{p}{\sqrt{1+p^2}},\frac{1}{\sqrt{1+p^2}})$.
It then implies from (\ref{S-3}) and  (\ref{3C:S-3}) that for any $q\not =0,$ $f(p)$ attains its maximal value, if and only if $p\in\R$ satisfies $\theta_p=0$ or $\pi$, where
\begin{equation}\label{3B:S-3}
	\begin{split}
\max_{p\in\R}f(p)&=\frac{{(q^2 + 1)^2(\mu_2 - \mu_3)(\mu_1 - \mu_3)}}{{\mu_2 - \mu_3 + q^2(\mu_1 - \mu_3)}}
 - \mu_1- \mu_2 q^2\\
&\quad-\frac{q^2(\mu_1-\mu_2)^2(\mu_1^2q^2-\mu_1\mu_3q^2+\mu_2^2-\mu_2\mu_3)}{[(\mu_3-\mu_2)+q^2(\mu_3-\mu_1)][(\mu_3-\mu_2)^2+(\mu_3-\mu_1)^2q^2]}\\
&=\frac{{(q^2 + 1)^2(\mu_2 - \mu_3)(\mu_1 - \mu_3)}}{{\mu_2 - \mu_3 + q^2(\mu_1 - \mu_3)}}
- \mu_2 q^2 - \mu_1\\
&\quad-\frac{q^2(\mu_1-\mu_2)^2}{(\mu_3-\mu_2)+q^2(\mu_3-\mu_1)}+
\frac{q^2\mu_3(\mu_1-\mu_2)^2}{(\mu_3-\mu_2)^2+q^2(\mu_3-\mu_1)^2}\\
&=-\mu_3\frac{[q^2(\mu_1-\mu_3)+(\mu_2-\mu_3)]^2}{(\mu_3-\mu_1)^2q^2+(\mu_3-\mu_2)^2}>0.
	\end{split}	
\end{equation}
Moreover, we deduce from (\ref{S-3}) and  (\ref{3C:S-3}) that for any $q\not =0,$ $f(p)$ is monotonically decreasing in $\theta _p\in [0,\pi/2)$,  and $f(p)$ is monotonically increasing in $\theta _p\in (\pi/2, \pi]$, where
\begin{equation}\label{3A:S-3}
\lim_{\theta_p\to \frac{\pi}{2}} f(p)=-\infty .
\end{equation}
We thus conclude from (\ref{3B:S-3}) and (\ref{3A:S-3}) that for any  $q\not= 0$, $f(p)$ admits exactly two different zero points denoted by $\underline{p}=\underline{p}(q)<\overline{p}=\overline{p}(q)$, $i.e.$,
\begin{equation}\label{SS-9}
	f(\underline{p})=f(\overline{p})=0,  \ \  \hbox{where}\ \  \underline{p}<\overline{p}.
\end{equation}
This further shows from (\ref{S-3}) that $0< u_3'(0)^2$, $i.e.$, $f(p)>0$, if and only if
\begin{equation}\label{N-3}
p\in(\underline{p},\,\,\overline{p}),\ \ \mbox{where $\ \underline{p}=\underline{p}(q)<\overline{p}=\overline{p}(q) \  $ satisfy $f(\underline{p})=f(\overline{p})=0$.}
\end{equation}
Therefore, we conclude that the system \eqref{I-1}--\eqref{I-3} at $x=0$
admits a solution
$$\big(u_1(0), u_2(0), u_3(0), u'_1(0), u'_2(0), u'_3(0)\big),$$
if and only if
$q$ and $p$  satisfy \eqref{N-17}.

Finally, since $u_3(0)=0$, if $q$ and $p$ satisfy \eqref{N-17}, then one can verify from (\ref{N-4}) and (\ref{N-1}) that \eqref{I-1}--\eqref{I-3} admits at most four solutions of the form (\ref{N-21}), due to the  symmetry of \eqref{N-21A} with respect to $u_1(0),u_2(0),u'_1(0),$ and $u'_3(0)$. This completes the proof of Proposition \ref{prop2B}.
\qed


For $\mu_1<\mu_2<\mu_3<0$, Proposition \ref{prop2B} implies that for $\tilde{u}_3(0)=0$ and $\frac{\tilde{u}_2(0)}{\tilde{u}_1(0)}=q\neq 0$, there exist {\em at most} four solutions $$\big(u_1(0),u_2(0),u_3(0),u'_1(0),u'_2(0),u'_3(0)\big)$$
of  \eqref{I-1}--\eqref{I-3} satisfying \eqref{N-4} and \eqref{N-17}. Since any solution $u=(u_1,u_2,u_3)$  of  \eqref{GPs} is uniquely determined by its initial condition
$\big(u_1(0),u_2(0),u_3(0),u'_1(0),u'_2(0),u'_3(0)\big),$
it follows from Proposition \ref{prop2B} that \eqref{GPs} admits at most four solutions satisfying \eqref{N-4} and \eqref{N-17}. For $\mu_1<\mu_2<\mu_3<0$, in the following we further prove that \eqref{GPs} admits exactly four solutions, whose initial conditions satisfy \eqref{N-4} and \eqref{N-17}. For this purpose, since any solution $\tilde{u}=(\tilde{u}_1,\tilde{u}_2,\tilde{u}_3)$ of \eqref{GPs} can be written as $\tilde{u}_i=\frac{g_i}{f}$ satisfying Lemma \ref{lem-1}, which is determined by $a_1,a_2,a_3$, it suffices to prove that for $\tilde{u}_3(0)=0$ and $\frac{\tilde{u}_2(0)}{\tilde{u}_1(0)}=q\neq 0$, there exist  {\em at least} four vectors $(a_1,a_2,a_3)$, where $a_i\neq 0$ holds for $i=1, 2, 3,$ so that
\begin{equation}\label{N-77}
	\big(\tilde{u}_1(0), \tilde{u}_2(0), 0, \tilde{u}'_1(0), \tilde{u}'_2(0), \tilde{u}'_3(0)\big) \ \, \mbox{satisfies \eqref{N-4} and \eqref{N-17}}.
\end{equation}

In order to prove   (\ref{N-77}), suppose $\tilde{u}=(\tilde{u}_1,\tilde{u}_2,\tilde{u}_3)$ is a solution of \eqref{GPs} and satisfies $\tilde{u}'_3(0)\neq 0$, $\tilde{u}_3(0)=0$ and $\frac{\tilde{u}_2(0)}{\tilde{u}_1(0)}=q\neq 0$. Applying Lemma \ref{lem-1}, $u$ can be written as $\tilde{u}_i=\frac{g_i}{f}$ in terms of  $a_1,a_2$ and $a_3$, and it then follows from \eqref{1.1} that $a_1,a_2$ and $a_3$ satisfy
\begin{equation}\label{N-8}
1+\frac{1}{4}\left(\frac{a_1^2(\eta_3-\eta_1)}{\eta_1^2(\eta_3+\eta_1)}+\frac{a_2^2(\eta_3-\eta_2)}{\eta_2^2(\eta_3+\eta_2)}\right)+\frac{1}{16}\frac{a_1^2a_2^2(\eta_1-\eta_2)^2(\eta_3-\eta_1)(\eta_3-\eta_2)}{\eta_1^2\eta_2^2(\eta_1+\eta_2)^2(\eta_3+\eta_1)(\eta_3+\eta_2)}=0,
\end{equation}
and
\begin{equation}\label{N-9}
\frac{a_2\left[1+\frac{1}{4}\Big(\frac{a_1^2(\eta_2-\eta_1)}{\eta_1^2(\eta_2+\eta_1)}+
\frac{a_3^2(\eta_2-\eta_3)}{\eta_3^2(\eta_2+\eta_3)}\Big)+\frac{1}{16}
\frac{a_1^2a_3^2(\eta_1-\eta_3)^2(\eta_2-\eta_1)(\eta_2-\eta_3)}{\eta_1^2\eta_3^2
(\eta_1+\eta_3)^2(\eta_2+\eta_1)(\eta_2+\eta_3)}\right]}{a_1\left[1+\frac{1}{4}
\Big(\frac{a_2^2(\eta_1-\eta_2)}{\eta_2^2(\eta_1+\eta_2)}+\frac{a_3^2(\eta_1-\eta_3)}{\eta_3^2(\eta_1+\eta_3)}\Big)
+\frac{1}{16}\frac{a_2^2a_3^2(\eta_2-\eta_3)^2(\eta_1-\eta_2)(\eta_1-\eta_3)}{\eta_2^2\eta_3^2(\eta_2+\eta_3)^2
(\eta_1+\eta_2)(\eta_1+\eta_3)}\right]}=q\neq 0,
\end{equation}
where and below $\eta_i=\sqrt{|\mu_i|}>0$ holds for $i=1, 2, 3$.
It yields from (\ref{N-9}) that $a_i\neq 0$ holds for all $i=1, 2, 3,$ due to the assumption $\tilde{u}_3'(0)\neq 0$.
Denote
\begin{equation}\label{N-18}
X=\frac{a_1}{\eta_1},\,\,Y=\frac{a_2}{\eta_2},\,\, Z=\frac{a_3}{\eta_3},
\end{equation}
and
\begin{equation}\label{N-19}
0<A_{12}=\frac{\eta_1-\eta_2}{\eta_1+\eta_2}<1,\,\,0<A_{13}=\frac{\eta_1-\eta_3}{\eta_1+\eta_3}<1,
\,\,0<A_{23}=\frac{\eta_2-\eta_3}{\eta_2+\eta_3}<1.
\end{equation}
Thus, the problem $(a_1,a_2,a_3)$ of \eqref{N-8} and \eqref{N-9} can be simplified as the problem $(X,Y,Z)$ of
\begin{equation}\label{N-14}
	1-\frac{1}{4}\big(A_{13}X^2+A_{23}Y^2\big)+\frac{1}{16}A_{12}^2A_{13}A_{23}X^2Y^2=0,
\end{equation}
and
\begin{equation}\label{N-16} \frac{\eta_2Y\left[1+\frac{1}{4}(-A_{12}X^2+A_{23}Z^2)-\frac{1}{16}A^2_{13}A_{12}
A_{23}X^2Z^2\right]}{\eta_1X\left[1+\frac{1}{4}(A_{12}Y^2+A_{13}Z^2)+\frac{1}{16}A^2_{23}A_{12}A_{13}Y^2Z^2\right]}=q,
\end{equation}
where $q\not =0$, $X\not =0$, $Y\not =0$, and $Z\not =0$.
Denote the intersection of \eqref{N-14} and \eqref{N-16} by
\begin{equation}\label{sp-s}
\begin{split}
S:=\big\{(X,Y,Z):\ \, & (X,Y,Z)\,\ \hbox{satisfies \eqref{N-14} and \eqref{N-16}}, \ \, \\
&\hbox{where}\,\  X\neq 0,\,\ Y\neq 0,\,\  \hbox{and}\ \, Z\neq 0\big\}.
\end{split}
\end{equation}
\textcolor{black}{It can be seen from Figure 2  that the ratio $\frac{a_2}{a_1}$ (or equivalently $\frac{Y}{X}$) can take any real value.}
We then obtain from (\ref{N-16}) that for any $q\neq 0$, there always exists a solution $(X=X(a_1),Y=Y(a_2),Z=Z(a_3))$ satisfying the above system \eqref{N-14} and \eqref{N-16}. Therefore, for any $q\neq 0$, we have
$S\neq \emptyset$.

\begin{figure}[t]
	\centering
	\includegraphics[width=0.9\textwidth]{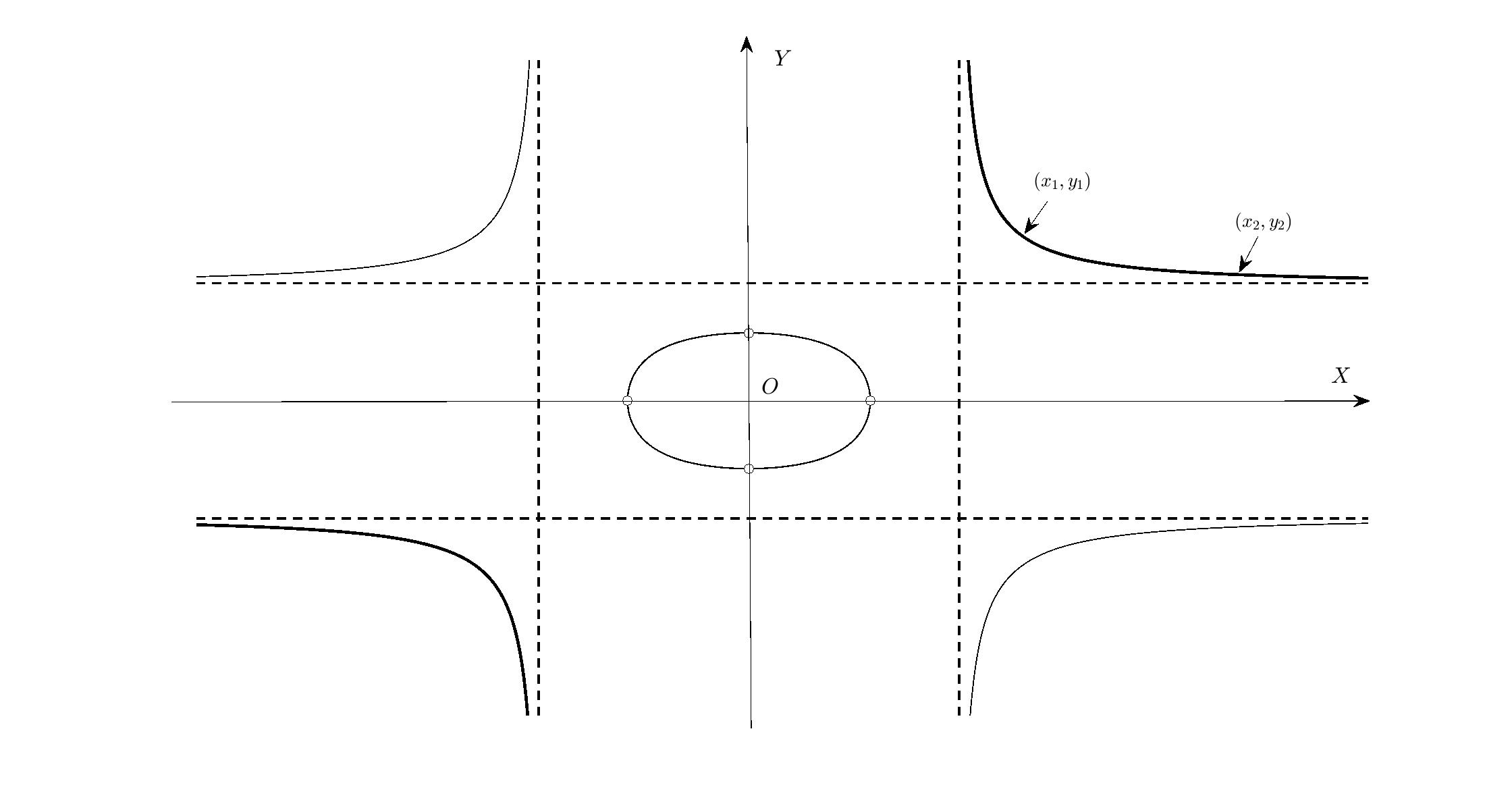}
	\caption{The curves $(X, Y)$ of \eqref{N-14}.}
	\label{fig-1} %
\end{figure}

For any fixed $q\neq 0$, we define
\begin{equation}\label{N-11}
p:=\frac{\tilde{u}'_2(0)}{\tilde{u}'_1(0)}, \ \ \mbox{so that (\ref{N-4}) is fully satisfied. }
\end{equation}	
In order to analyze $p$, it is obvious from Lemma \ref{lem-1} that
\[	
p\ \ \hbox{is a continuous function of $\big(X(a_1),Y(a_2), Z(a_3)\big)$.}\]
Moreover, for any fixed $q\neq 0$, $p(X,Y,Z)$ can be thought of as:
\begin{equation}
p: \ \, \hbox{the curve $S$}\to  \R^1\ \ \hbox{is continuous.}
\end{equation}
The following lemma shows that for any fixed $q\neq 0$, the range of $p(X,Y,Z)$ along at least four connected branches of the curve $S$ defined by \eqref{sp-s} is equal to $(\underline{p},\,\,\overline{p})$, where $\underline{p}=\underline{p}(q)<\overline{p}=\overline{p}(q)$ are as in (\ref{S-9}). This implies equivalently that for $\tilde{u}_3(0)=0$ and $\frac{\tilde{u}_2(0)}{\tilde{u}_1(0)}=q\neq 0$, there exist  {\em at least} four vectors $(a_1,a_2,a_3)$, where $a_i\neq 0$ holds for all $i=1, 2, 3,$ so that $ \big(\tilde{u}_1(0), \tilde{u}_2(0), 0, \tilde{u}'_1(0), \tilde{u}'_2(0), \tilde{u}'_3(0)\big)$ also satisfies \eqref{N-17}. Once Lemma \ref{lemS} below is proved,  this therefore completes the proof of (\ref{N-77}) in view of (\ref{N-11}).


\begin{lem}\label{lemS} For any fixed  $q\neq 0$,
let the curve $S$ be defined by \eqref{sp-s}. Then $S$ has four connected branches $S_1$, $S_2$, $S_3$ and $S_4$, i.e., $S_1\cup S_2\cup S_3\cup S_4\subset S$, along each of which the function $p$ defined in \eqref{N-11} is onto $(\underline{p},\,\,\overline{p})$, i.e., $p(S_i)=(\underline{p},\,\,\overline{p})$ holds for $i=1,2,3,4$.
\end{lem}

\noindent{\bf Proof.} In view of (\ref{N-18}) and (\ref{N-19}), we rewrite \eqref{N-14}  as
\begin{equation}\label{SSP-2}
\Big(\frac{1}{4}A_{13}A_{12}X^2-\frac{1}{A_{12}}\Big)
\Big(\frac{1}{4}A_{23}A_{12}Y^2-\frac{1}{A_{12}}\Big)=\frac{1}{A^2_{12}}-1>0,\ \, X\neq 0, \ \, Y\neq 0,
\end{equation}
which has eight connected branches in the XY plane, see above Figure 2.
On the other hand,  \eqref{N-16} can be regarded as
\begin{equation}\label{sp-3}
\begin{aligned}
Z^2:&=F(X,Y)\\
&=\frac{-4 \left( A_{12} X Y^{2} \eta_{1} q +  A_{12} X^{2} Y \eta_{2} + 4X \eta_{1} q - 4 Y \eta_{2}\right)}{  \left(  A_{12} A_{13} A_{23}^{2} X Y^{2} \eta_{1} q +   A_{12} A_{13}^{2} A_{23} X^{2} Y \eta_{2} + 4 A_{13} X \eta_{1} q - 4 A_{23} Y \eta_{2}\right)},
\end{aligned}
\end{equation}
where $F(X,Y)$ is a rational function of $X$ and $Y$, and continuous with respect to $X$ and $Y$ along each connected branch of \eqref{N-14}.

Since $S\neq\emptyset$,  we consider from (\ref{SSP-2}) that the point $(x_0,y_0,z_0)\in S$ satisfies
\begin{equation}\label{SSP-1}
z_0>0,\,\ \frac{1}{4}A_{13}A_{12}x_0^2-\frac{1}{A_{12}}>0\,\ \hbox{and}\,\ \frac{1}{4}A_{23}A_{12}y_0^2-\frac{1}{A_{12}}>0,
\end{equation}
which corresponds to a point $(x_0,y_0)$ on the bold black curve in the first quadrant of Figure 2. Set
\begin{equation}\label{sp-1}
S_1:=\big\{(x_0,y_0,z_0)\big\}\cup\big\{(x,y,z):\,\ (x,y,z)\in S\,\ \hbox{is connected with}\,\ (x_0,y_0,z_0)\big\},
\end{equation}
where $(x_0,y_0,z_0)$ is as in (\ref{SSP-1}). We now prove that
\begin{equation}\label{sp-2}
\sup\big\{x^2+y^2:\,\ (x,y,z)\in S_1\subset S\big\}<\infty.
\end{equation}
On the contrary, suppose \eqref{sp-2} is false, then there exists a sequence $\{(x_n,y_n,z_n)\}\subset S_1$ such that either $\lim_{n\to\infty}x_n^2=\infty$ or $\lim_{n\to\infty}y_n^2=\infty$.
Without loss of generality, we assume that $\lim_{n\to\infty}x_n^2=\infty$. Since $(x_n,y_n)$ satisfies \eqref{N-14}, we obtain from (\ref{SSP-2}) that the sequence $\{y_n\}$ is bounded uniformly in $n$. Thus, we derive from \eqref{sp-3} that
\[
\lim_{n\to\infty}z_n^2=-\frac{4}{A^2_{13}A_{23}}<0,
\]
which is impossible. Thus, \eqref{sp-2} holds true.

Define the left endpoint $T_1=(x_1,y_1,z_1)$ and the right endpoint $T_2=(x_2,y_2,z_2)$ of $S_1$ as follows:
\begin{equation}\label{sp-4}
\begin{split}
x_1=\inf\big\{x:\ (x,y,z)\in S_1\big\},\,\ &y_1=\sup\big\{y:\ (x,y,z)\in S_1\big\},\,\    z_1=\sqrt{F(x_1,y_1)},\\
x_2=\sup\big\{x:\ (x,y,z)\in S_1\big\},\,\ &y_2=\inf\big\{y:\ (x,y,z)\in S_1\big\}, \,\  z_2=\sqrt{F(x_2,y_2)},
\end{split}
\end{equation}
where $z_i=+\infty$ is defined for $i=1,2$, provided that $(x_i,y_i)$ is a singular point of $F(X,Y)$. We next claim that the point  $z_i$ of (\ref{sp-4}) satisfies
\begin{equation}\label{sp-6}
z_i\in\{0,+\infty\},\,\ \hbox{where}\,\ i=1,2.
\end{equation}
Indeed, on the contrary, suppose \eqref{sp-6} is false. Without loss of generality, then we may assume from (\ref{SSP-1}) that $z_1\in(0,+\infty)$ and $F(x_1,y_1)=z_1^2$. Taking a small enough $\eps>0$ so that $z^2_1-\eps>0$. Since $F(X,Y)$ is continuous with respect to the curve  \eqref{N-14}, we obtain that
\[\begin{split}
&F^{-1}\big((z^2_1-\eps,z^2_1+\eps)\big)\\
:=&\big\{(x,y)\,\ \hbox{satisfies}\, \eqref{N-14} \,\ \hbox{and}\, (\ref{SSP-1}):\ F(x,y)\in (z^2_1-\eps,z^2_1+\eps)\big\}
\end{split}
\]
is a one dimensional open set of the curve  \eqref{N-14}. Moreover, for any $(x_3,y_3)\in F^{-1}\big((z^2_1-\eps,z^2_1+\eps)\big)$, it is obvious that
$\big (x_3, y_3, \sqrt{F(x_3,y_3)}\big)\in S_1$. Since $(x_1,y_1)\in F^{-1}((z^2_1-\eps,z^2_1+\eps))$,  $(x_1,y_1,z_1)$ is an interior point of $S_1$, which however contradicts with the fact that $(x_1,y_1,z_1)$ is the left endpoint of $S_1$. Therefore, the claim \eqref{sp-6} holds true.

It yields from \eqref{sp-6} that if
$$(X,Y,Z):=\Big(\frac{a_1}{\eta_1},\frac{a_2}{\eta_2},\frac{a_3}{\eta_3}\Big)=(x_i,y_i,z_i),\,\ i=1,2,$$
where the point $(x_i,y_i,z_i)$ is as in (\ref{sp-4}), then
the expression of $\tilde u'_3$ in Lemma \ref{lem-1}  gives that $\tilde u_3'(0)=0$, which implies from \eqref{S-3} that $f(p)=\frac{\tilde u_3'(0)^2}{\tilde u_1(0)^2}=0$. This further yields from \eqref{S-9} that the point $(x_i,y_i,z_i)$ of (\ref{sp-4}) satisfies either $p(x_i,y_i,z_i)=\overline{p}$ or $p(x_i,y_i,z_i)=\underline{p}$ for $i=1,2$, where $\overline{p}<\underline{p}$ are as in \eqref{S-9}.

We now claim that
\begin{equation}\label{S-1}
	\big\{p(x_1,y_1,z_1),\,\,p(x_2,y_2,z_2)\big\}=\big\{\underline{p},\,\,\overline{p}\big\},
\end{equation}
where the point $(x_i,y_i,z_i)$ is as in (\ref{sp-4}) for $i=1, 2$, and $\overline{p}<\underline{p}$ are as in \eqref{S-9}.
Indeed, suppose \eqref{S-1} is false. Without loss of generality, then we may assume from (\ref{S-1}) that $p(x_1,y_1,z_1)=p(x_2,y_2,z_2)=\underline{p}$. Let $M=\max_{(x,y,z)\in S_1}p(S_1)>\underline{p}$. For any $p_0\in(\underline{p},M)$, we then obtain that there exist at least two points $(x_4,y_4,z_4)\in S_1$ and $(x_5,y_5,z_5)\in S_1$ such that $p(x_4,y_4,z_4)=p(x_5,y_5,z_5)=p_0$.
Due to the symmetry of $p$ in \eqref{N-11} with respect to $z$, we obtain that $p(x_4,y_4,-z_4)=p(x_5,y_5,-z_5)=p_0$.
On the other hand, the expressions of $\tilde u_1(x)$  and $\tilde u_2(x)$  in  Lemma \ref{lem-1} yield that $p(x,y,z)=p(-x,-y,z)$. By the symmetry of \eqref{N-14} and \eqref{N-16} with respect to $(X,Y,Z)$, we thus obtain that $(-x_4,-y_4,\pm z_4)$, $(-x_5,-y_5,\pm z_5)$, $(x_4,y_4,-z_4)$, $(x_5,y_5,-z_5)\in S$, and
$$p(-x_4,-y_4,\pm z_4)=p(-x_5,-y_5,\pm z_5)=p(x_4,y_4,-z_4)=p(x_5,y_5,-z_5)=p_0.$$
We hence conclude from above that $p(x,y,x)=p_0$ has at least eight solutions, which however contradicts with Proposition \ref{prop2B}. Therefore, the claim \eqref{S-1} holds true.

Applying \eqref{S-1}, we conclude immediately that  $p(S_1)=(\underline{p},\,\,\overline{p})$, due to the continuity  of $p$ with respect to $S_1$. By the symmetry of $S$ and  $p(S_1)=(\underline{p}, \overline{p})$, we obtain from   Proposition \ref{prop2B} that $S$ has exactly four connected branches, where $S_1$ is given by \eqref{sp-1},
\[
S_2:=\big\{(x,y,z): \ (x,y,-z)\in S_1\big\},\,\ S_3:=\big\{(x,y,z): \ (-x,-y,-z)\in S_1\big\},
\]
and
\[
S_4:=\big\{(x,y,z): \ (-x,-y,z)\in S_1\big\}.
\]
Here the curve $S_1$ of the  $XYZ$ space satisfies $Z>0$, and its projection onto the $XY$ plane is the bold black curve in the first quadrant of Figure 2, where the left endpoint is  at
$(x_1,y_1)$ and the right endpoint is at
$(x_2,y_2)$. $S_2$ is  the reflection of $S_1$ about the $XY$ plane. $S_3$ and $S_4$ are the symmetric curves of
$S_1$ and $S_2$ with respect to the origin, respectively.
Moreover, we have $p(S_i)=p(S_1)=(\underline{p}, \overline{p})$ for $i=1,2,3,4$, which completes the proof of Lemma \ref{lemS}.     \qed


\subsection{Proofs of main results}

In this subsection, we first establish  Theorem \ref{thm1.1} on the classification of general solutions for  (\ref{GPs}). As the applications of Theorem \ref{thm1.1}, we then address the proofs of Theorem \ref{thm1.2} and  Corollary \ref{thm1.3}.

\vspace{.1cm}

\noindent\textbf{Proof of Theorem \ref{thm1.1}.} Suppose $u=(u_1, u_2, u_3)$ is a solution of (\ref{GPs}). If either $\mu_1=\mu_2<\mu_3<0$ or $\mu_1<\mu_2=\mu_3<0$ or $\mu_1=\mu_2=\mu_3<0$, then we derive from Proposition \ref{prop2A} that Theorem \ref{thm1.1} holds true. Thus, the rest is to consider $\mu_1<\mu_2<\mu_3<0$, which implies that $u_3$ changes sign. \textcolor{black}{Indeed, if either $u_3>0$ or $u_3<0$, then $\mu_3$ is the first eigenvalue of the operator $-\frac{d^2}{dx^2}-2(u_1^2+u_2^2+u_3^2)\ $ in $L^2(\R)$, which implies that $\mu_3\leq\mu_2\leq\mu_1$, a contradiction. Hence, $u_3$ vanish somewhere, and up to a translation, we may assume that $u_3(0)=0$.} Moreover, we may assume that $ u_3'(0)\neq 0$, since otherwise we have $u_3\equiv 0$ and hence Theorem \ref{thm1.1} follows directly from \cite[Lemma 15]{Lewin}. Therefore, if $\mu_1<\mu_2<\mu_3<0$, then the proof of Theorem \ref{thm1.1} is divided into the following four different cases:
\begin{itemize}
\item[ (i).]  $u_i(0)=u_3(0)=0$ holds for $i=1$ or $2$.

\item[(ii).]  $u'_1(0)=u'_2(0)=u_3(0)=0$.

\item[(iii).]  $ \frac{u_2(0)}{u_1(0)}=q\neq 0,\,\,  u'_2(0)=pu'_1(0),\ \  u_1'(0)\neq 0,\,\ u_3(0)=0,\,\,  u_3'(0)\neq 0$.

\item[(iv).]
$ \frac{u_1(0)}{u_2(0)}=q\neq 0,\,\,  u'_1(0)=pu'_2(0),\ \  u_2'(0)\neq 0,\,\  u_3(0)=0,\,\, u_3'(0)\neq 0$.
\end{itemize}
In the first two cases (i) and (ii), we derive from Lemma \ref{lem:3.1} that either $u_1\equiv 0$ or $u_2\equiv 0$. Therefore,  a similar argument of \cite[Lemma 15]{Lewin} then yields Theorem \ref{thm1.1}, and we are done. Moreover, since the proof of case (iv) is very similar to that of case (iii),  the rest is to address the proof of case (iii).

Indeed, suppose that case (iii) happens. It then follows from
Proposition \ref{prop2B} that for any $q\neq 0$ and $p\in(\underline{p},\,\,\overline{p})$, there exist  at most four solutions of \eqref{GPs}. On the other hand, we obtain from Lemma \ref{lemS} that for any $q\neq 0$ and $p\in(\underline{p},\,\,\overline{p})$, there exist at least four solutions of \eqref{GPs} satisfying \eqref{exp1}. We thus conclude from above that for any  $q\neq 0$ and $p\in(\underline{p},\,\,\overline{p})$, there exist exactly four solutions for \eqref{GPs}, which must satisfy \eqref{exp1}. This completes the proof of Theorem \ref{thm1.1}.    \qed

To discuss the applications of Theorem \ref{thm1.1}, we next give the following mass analysis of solutions for \eqref{GPs}.

\begin{lem}\label{lem-2}
For $\mu_1<\mu_2<\mu_3<0$, suppose $(u_1,u_2,u_3)\in H^1(\R)^3$ is a  solution of \eqref{GPs} satisfying $u_i\not\equiv 0$ for all $i=1,2,3$. Then it must have $\int_{\R}u^2_i(x)dx=2\eta_i$ for $i=1,2,3$.
\end{lem}

\noindent{\bf Proof.} Following Theorem \ref{thm1.1}, we claim that the identity
\begin{equation} \label{N.4}
\begin{split}
&u^2_i(x)\\
=& -2\eta_i\Big(\frac{1+\sum_{j=1,j\neq i}^{3}\frac{ a_j^2}{4\eta^2_j}e^{2\eta_jx}+\sum_{j\neq i,k\neq i,1\leq j< k\leq 3}\frac{a_j^2a_k^2(\eta_j-\eta_k)^2}{16\eta_j^2\eta_k^2(\eta_j+\eta_k)^2}e^{2(\eta_j+\eta_k)x}}{f(x)}\Big)'
\end{split}
\end{equation}
holds true for all $i=1,2,3$, where $f(x)$ is as in Theorem \ref{thm1.1}. One can check that if the identity (\ref{N.4}) holds for all $i=1,2,3$,  then we have
\[
\int_{\R}u^2_i(x)dx=2\eta_i,\ \ i=1, 2, 3,
\]
which thus proves the lemma.

In the following, we only need to prove the claim (\ref{N.4}). Without loss of generality, we just prove  the claim (\ref{N.4}) for $i=1$.
Actually, we note from Theorem \ref{thm1.1} that
\[
\begin{aligned} g_1(x)&=a_1e^{\eta_1x}\Big\{1+\frac{a_2^2(\eta_1-\eta_2)}{4\eta^2_2(\eta_1+\eta_2)}e^{2\eta_2x}+
\frac{a_3^2(\eta_1-\eta_3)}{4\eta^2_3(\eta_1+\eta_3)}e^{2\eta_3x}\\ &\qquad\qquad \,+\frac{a_2^2a_3^2(\eta_1-\eta_2)(\eta_1-\eta_3)(\eta_2-\eta_3)^2}{16\eta_2^2\eta_3^2(\eta_1+\eta_2)
(\eta_1+\eta_3)(\eta_2+\eta_3)^2}e^{2(\eta_2+\eta_3)x}\Big\}:=a_1e^{\eta_1x}\hat{g}(x).
\end{aligned}
\]
Define
\[
\underline{f}(x)=1+\frac{a_2^2}{4\eta^2_2}e^{2\eta_2x}+\frac{a_3^2}{4\eta^2_3}e^{2\eta_3x}+\frac{a_2^2a_3^2
(\eta_2-\eta_3)^2}{16\eta_2^2\eta_3^2(\eta_2+\eta_3)^2}e^{2(\eta_2+\eta_3)x},
\]
and
\[
\begin{aligned}
\overline{f}(x)&=f(x)-\underline{f}(x)\\
&=\frac{a_1^2}{4\eta^2_1}e^{2\eta_1x}\Big\{1+\frac{a_2^2(\eta_1-\eta_2)^2}
{4\eta^2_2(\eta_1+\eta_2)^2}e^{2\eta_2x}+\frac{a_3^2(\eta_1-\eta_3)^2}{4\eta^2_3(\eta_1+\eta_3)^2}e^{2\eta_3x}\\
&\quad\qquad\qquad\ +\frac{a_2^2a_3^2(\eta_1-\eta_2)^2(\eta_1-\eta_3)^2(\eta_2-\eta_3)^2}{16\eta_2^2\eta_3^2(\eta_1+\eta_2)^2
(\eta_1+\eta_3)^2(\eta_2+\eta_3)^2}e^{2(\eta_2+\eta_3)x}\Big\}.
\end{aligned}
\]
We also set
\begin{equation}\label{N.5}
[F(x),G(x)]:=F'(x)G(x)-F(x)G'(x),
\end{equation}
so that
\begin{equation}\label{sp-7}
\begin{split}
[F(x),F(x)]=0,\,\ [F(x),G(x)]=-[G(x),F(x)],\\
 \hbox{and}\ \  [F(x),H(x)F(x)]=-H'(x)F^2(x).\qquad
 \end{split}
\end{equation}
We then derive from above that
\[
\mbox{RHS of}\ \, (\ref{N.4})=-\frac{2\eta_1}{f^2}\big[\,\underline{f}(x),\overline{f}(x)\big],
\]
and
\[
\mbox{LHS of}\ \, (\ref{N.4})=-\frac{2\eta_1}{f^2}\Big[\hat{g}(x), \frac{a_1^2}{4\eta_1^2}e^{2\eta_1x}\hat{g}(x)\Big].
\]
Applying \eqref{sp-7}, one can calculate that
\[
\big[\,\underline{f}(x),\overline{f}(x)\big]=\Big[\hat{g}(x), \frac{a_1^2}{4\eta_1^2}e^{2\eta_1x}\hat{g}(x)\Big].
\]
This proves  the claim \eqref{N.4}, and it therefore completes the proof of Lemma \ref{lem-2}.
\qed

By employing Theorem \ref{thm1.1}, we are now ready to address the proofs of Theorem  \ref{thm1.2} and  Corollary \ref{thm1.3}.

\vspace{.1cm}

\noindent\textbf{Proof of Theorem \ref{thm1.2}.}  Firstly, if $\mu_1=\mu_2=\mu_3=-\frac{9}{4}$, then we derive from Proposition \ref{prop2A} and Theorem \ref{thm1.1}  that
\eqref{GPs} admits normalized solutions $u=(u_1,u_2,u_3)$ satisfying (\ref{exp1}), where  $|a_1|=|a_2|=|a_3|=1$,  and $\eta_1=\eta_2=\eta_3=\frac{3}{2}$. Thus, up to the translation, $u=(u_1,u_2,u_3)$ must be unique, in the sense that
\begin{equation*}
	u_1(x)\equiv\pm\frac{\sqrt{3}}{2\cosh(\frac{3x}{2})},\,\, u_2(x)\equiv\pm\frac{\sqrt{3}}{2\cosh(\frac{3x}{2})},\,\,u_3(x)\equiv\pm\frac{\sqrt{3}}{2\cosh(\frac{3x}{2})}
\end{equation*}
hold for three uncorrelated signs $\pm$.

Similarly,  if $\mu_1=\mu_2<\mu_3<0$, then we derive from Proposition \ref{prop2A} that any solution $u=(u_1,u_2,u_3)$ of \eqref{GPs} satisfies $u_1\equiv ku_2$,  $\int_{\R}(u^2_1+u^2_2)dx=2\eta_1$ and $\int_{\R}(u^2_3)dx=2\eta_3$. Further, if $\int_{\R}u^2_1dx=\int_{\R}u^2_2dx=\int_{\R}u^2_3dx=1$, then we have $u_1\equiv\pm u_2$, $\mu_1=\mu_2=-1$ and $\mu_3=-\frac{1}{4}$. In this case, we thus deduce from Theorem \ref{thm1.1} that \eqref{GPs} admits infinitely many normalized solutions $u=(u_1,u_2,u_3)$, which satisfy  (\ref{thm1:M8}) for all $A\neq 0$ and $B\neq 0$.

Additionally, if $\mu_1<\mu_2=\mu_3<0$, then we obtain from Proposition \ref{prop2A} that $\int_{\R}u^2_1dx=2\eta_1>\int_{\R}(u^2_2+u^2_3)dx=2\eta_2$, and hence $\int_{\R}u^2_1dx=\int_{\R}u^2_2dx=\int_{\R}u^2_3dx=1$ cannot occur, which implies the nonexistence of normalized solutions for \eqref{GPs}. Finally, if $\mu_1<\mu_2<\mu_3<0$, then we deduce from   Theorem \ref{thm1.1} and Lemma \ref{lem-2} that $\int_{\R}u^2_1dx=\int_{\R}u^2_2dx=\int_{\R}u^2_3dx=1$ cannot occur yet, which also implies the nonexistence of normalized solutions for \eqref{GPs}. This completes the proof of Theorem \ref{thm1.2}. \qed

\vspace{.1cm}

\noindent\textbf{Proof of Corollary \ref{thm1.3}.} (i). Suppose $u=(u_1,u_2,u_3)\in H^1(\R)^3$ is a least energy normalized solution of \eqref{GPs} \textcolor{black}{with fixed constants $\mu_1,\mu_2,\mu_3$}, and let $E(\mu_1,\mu_2,\mu_3,u_1,u_2,u_3)$ be the corresponding energy of the solution  $u$, $i.e.$,
\begin{equation}\label{sp-11:3}
	E(\mu_1,\mu_2,\mu_3,u_1,u_2,u_3):=	\sum_{k=1}^3\int_{\R}(u_k')^2dx-\int_{\R}\Big(\sum_{k=1}^3(u_k)^2\Big)^2dx.
\end{equation}
\textcolor{black}{Since $(u_1,u_2,u_3)$ depends on $\mu_1,\mu_2$ and $\mu_3$ in view of \eqref{GPs}, $E(\mu_1,\mu_2,\mu_3,u_1,u_2,u_3)$ also depends on $\mu_1,\mu_2$ and $\mu_3$.}
Following Theorem \ref{thm1.2}, it suffices to show that
\begin{equation}\label{sp-011:3}
E(-\frac{9}{4},-\frac{9}{4},-\frac{9}{4},u_1,u_2,u_3)<E(-1,-1,-\frac{1}{4},u_1,u_2,u_3).
\end{equation}
Indeed, if (\ref{sp-011:3}) holds true, then it implies that Theorem \ref{thm1.2} (i) must occur. We thus conclude from Theorem \ref{thm1.2} (i) that it necessarily has $\mu_1=\mu_2=\mu_3=-\frac{9}{4}$, and $u$ must satisfy \eqref{2.0}.

To prove (\ref{sp-011:3}), we multiply  the $j$-th equation in \eqref{GPs} by $u_j$, sum up the resulting identities, and finally integrate it on $\R$. We then obtain that
\begin{equation}\label{sp-12:3}
-\sum_{k=1}^3\int_{\R}	(u_k')^2dx+2\int_{\R}\Big(\sum_{k=1}^3(u_k)^2\Big)^2dx=-(\mu_1+\mu_2+\mu_3).
\end{equation}
On the other hand, we derive from \eqref{I-1} that
\begin{equation}\label{sp-13:3}
\sum_{k=1}^3\int_{\R}	 (u_k')^2dx+\int_{\R}\Big(\sum_{k=1}^3u_k^2\Big)^2dx=-(\mu_1+\mu_2+\mu_3).
\end{equation}
It then yields from \eqref{sp-12:3} and \eqref{sp-13:3} that
\[
\sum_{k=1}^3\int_{\R} (u_k')^2dx=-\frac{1}{3}(\mu_1+\mu_2+\mu_3),\,\ \hbox{and}\,\ \int_{\R}\Big(\sum_{k=1}^3u_k^2\Big)^2=-\frac{2}{3}(\mu_1+\mu_2+\mu_3),
\]
which further implies that
\[
E(\mu_1,\mu_2,\mu_3,u_1,u_2,u_3)=\frac{1}{3}(\mu_1+\mu_2+\mu_3).
\]
Direct computations thus yield that $$E(-\frac{9}{4},-\frac{9}{4},-\frac{9}{4},u_1,u_2,u_3)=-\frac{9}{4}<E(-1,-1,-\frac{1}{4},u_1,u_2,u_3)=-\frac{3}{4},$$
which proves  (\ref{sp-011:3}), and we are  done.

(ii). Following Theorem \ref{thm1.2}, we only need to prove the nonexistence of orthonormal solutions for the following two cases:
\begin{itemize}
	\item[(i).] $\mu_1=\mu_2=\mu_3=-\frac{9}{4}$;
	\item[(ii).] $\mu_1=\mu_2=-1$ and $\mu_3=-\frac{1}{4}$.
\end{itemize}
Actually, one can directly conclude from Proposition \ref{prop2A} that $u_1$ and $u_2$ cannot be orthogonal in $L^2(\R)$ for each of above two cases, which hence proves Corollary \ref{thm1.3}. \qed

%
%
%
%
%

\section{Nondegeneracy of the System \eqref{GPs}}

The purpose of  this section is to establish Theorem \ref{thm-nd} on the nondegeneracy of solutions $u=(u_1,u_2,u_3)$ for the system \eqref{GPs} satisfying $u_i\not\equiv 0$ for $i=1,2,3$. Towards this aim, suppose that $u=(u_1,u_2,u_3)$ satisfies the system \eqref{GPs}, where $\mu_1<\mu_2<\mu_3<0$ are arbitrary. We recall from  \eqref{1:op} that
 \begin{equation*}
L_u:= -\frac{d^2}{dx^2}-2(u_1^2+u_2^2+u_3^2)\ \ \mbox{in}\,\ L^2(\R),
 \end{equation*}
and  consider the following linearized system of \eqref{GPs} around the solution $u$:
\begin{equation}\label{Line}
	\left\{\begin{array}{lll}
		\phi_1''+2(u_1^{2}+u_2^{2}+u_3^{2})\phi_1+4(u_1\phi_1+u_2\phi_2+u_3\phi_3)u_1=-\mu_1\phi_1 \ \ \mbox{in}\, \ \R,\\[3mm]
		\phi_2''+2(u_1^{2}+u_2^{2}+u_3^{2})\phi_2+4(u_1\phi_1+u_2\phi_2+u_3\phi_3)u_2=-\mu_2\phi_2 \ \ \mbox{in}\, \ \R,\\[3mm]
		\phi_3''+2(u_1^{2}+u_2^{2}+u_3^{2})\phi_3+4(u_1\phi_1+u_2\phi_2+u_3\phi_3)u_3=-\mu_3\phi_3 \ \ \mbox{in}\, \ \R.
	\end{array}\right.
\end{equation}

We first address the following lemma on the constants of motion for $(\phi_1,\phi_2,\phi_3)$ satisfying (\ref{Line}), which is essentially the linearized version of Lemma  \ref{lem2.2}.

\begin{lem}\label{lem5.1}
Suppose $\phi =(\phi_1,\phi_2,\phi_3)\in H^1(\R)^3$ is a  solution of \eqref{Line}, then we have the following three constants of motion:	
\begin{equation}\label{D-1}
	\begin{split}
		&u_1'\phi_1'+u_2'\phi_2'+u_3'\phi_3'+2(u_1^2+u_2^2+u_3^2)(u_1\phi_1+u_2\phi_2+u_3\phi_3)\\
		&+\mu_1u_1\phi_1+\mu_2u_2\phi_2+\mu_3u_3\phi_3=0 \ \ \mbox{in}\, \ \R,	
	\end{split}
\end{equation}
\begin{equation}\label{D-2}
	\begin{split}
		&\sum_{1\leq j<k\leq3}(u_j'u_k-u_ju_k')(\phi_j'u_k+u_j'\phi_k-\phi_ju_k'-u_j\phi_k')\\[2mm]
		&\,\,+(u_1^2+u_2^2+u_3^2)\big[(\mu_2+\mu_3)u_1\phi_1+(\mu_1+\mu_3)u_2\phi_2+(\mu_1+\mu_2)u_3\phi_3\big]\\[2mm]
		&+(u_1\phi_1+u_2\phi_2+u_3\phi_3)\big[(\mu_2+\mu_3)u_1^2+(\mu_1+\mu_3)u_2^2+(\mu_1+\mu_2)u_3^2\big]\\[2mm]
		&\,\,+(\mu_2+\mu_3)u'_1\phi_1'+(\mu_1+\mu_3)u'_2\phi_2'+(\mu_1+\mu_2)u'_3\phi_3'\\[2mm]
		&\,\,+\mu_1(\mu_2+\mu_3)u_1\phi_1+\mu_2(\mu_1+\mu_3)u_2\phi_2+\mu_3(\mu_1+\mu_2)u_3\phi_3=0 \ \ \mbox{in}\, \ \R,
	\end{split}
\end{equation}
and
\begin{equation}\label{D-3}
	\begin{split}
		& \sum_{1\leq j<k\leq3,\,i\neq j,\,i\neq k}\mu_i(u_j'u_k-u_ju_k')(\phi_j'u_k+u_j'\phi_k-\phi_ju_k'-u_j\phi_k')\\[2mm]
		&\,\,+(u_1^2+u_2^2+u_3^2)(\mu_2\mu_3u_1\phi_1+\mu_1\mu_3u_2\phi_2+\mu_1\mu_2u_3\phi_3)\\[2mm]
		&+(u_1\phi_1+u_2\phi_2+u_3\phi_3)(\mu_2\mu_3u_1^2+\mu_1\mu_3u_2^2+\mu_1\mu_2u_3^2)\\[2mm]
		 &\,\,+\mu_2\mu_3u'_1\phi_1'+\mu_1\mu_3u'_2\phi_2'+\mu_1\mu_2u'_3\phi_3'+\mu_1\mu_2\mu_3(u_1\phi_1+u_2\phi_2+u_3\phi_3)=0 \ \ \mbox{in}\, \ \R.
	\end{split}
\end{equation}
\end{lem}

\noindent{\bf Proof.} Since $\phi =(\phi_1,\phi_2,\phi_3)\in H^1(\R)^3$ is a  solution of \eqref{Line}, we deduce from \eqref{Line} that $\phi_j$ and $\phi_j'$ vanish at infinity, where $j=1,2,3$. We then get that $\phi_j$, $u_j$, $\phi_j'$ and $u_j'$ vanish at infinity, where $j=1,2,3$.

We now prove \eqref{D-1}--\eqref{D-3} as follows. \eqref{D-1} is derived  by using $u_j'$ and $\phi_j'$ in the following way:
\[
0=\inte \big[\text{(}\ref{Line}\text{$)_j$} \times u'_j + \text{(}\ref{GPs}\text{$)_j$} \times \phi_j'\big] dx.
\]
However, \eqref{D-2} is derived by choosing pairwise coupled integrating factors as follows:
\[
\begin{aligned}
0=&\inte \big[\text{(}\ref{Line}\text{$)_j$} \times\sum_{1\leq k\leq 3,\,k\neq j }\big[(u_j'u_k-u_ju_k')u_k+\mu_ku_j'\big] \\
&\quad + \text{(}\ref{GPs}\text{$)_j$} \times \sum_{1\leq k\leq 3,\,j\neq k}\big[(\phi_j'u_k-\phi_ju_k'+\phi_ku_j'-\phi_k'u_j)u_k+(u_j'u_k-u_ju_k')\phi_k+\mu_k\phi_j'\big]dx.
\end{aligned}
\]
Finally, \eqref{D-3} is obtained analogously by using integrating factors involving three components in the following way:
\[
\begin{aligned}
0=&\inte \big[\text{(}\ref{Line}\text{$)_j$} \times\sum_{1\leq i<k\leq 3,\,i\neq j,\,k\neq j}\big[\mu_i(u_j'u_k-u_ju_k')u_k+\mu_k(u_j'u_i-u_ju_i')u_i+\mu_i\mu_ku_j'\big]\big] \\
	&\quad + \text{(}\ref{GPs}\text{$)_j$} \times \sum_{1\leq i<k\leq 3,\,i\neq j,\,k\neq j}\big[\mu_i(\phi_j'u_k-\phi_ju_k'+\phi_ku_j'-\phi_k'u_j)u_k+\mu_i(u_j'u_k-u_ju_k')\phi_k\\
	&\quad+\mu_k(\phi_j'u_i-\phi_ju_i'+\phi_iu_j'-\phi_i'u_j)u_i+\mu_k(u_j'u_i-u_ju_i')\phi_i
	+\mu_i\mu_ku_j'\big].
\end{aligned}
\]
This proves Lemma \ref{lem5.1}. \qed

%
%
%
%
%
%
%
%
%
%


\begin{prop}\label{prop5.2} For $\mu_1<\mu_2<\mu_3<0$, let
\begin{equation}\label{D}
\begin{split}
\mathcal{D}:=&\big\{(\phi_1,\phi_2,\phi_3)\in H^1(\R)^3:\,
 (\phi_1,\phi_2,\phi_3)\ \, \hbox{is a
	solution of}\,\  (\ref{D-1})-(\ref{D-3})\big\},
\end{split}
\end{equation}
then $\mathcal{D}$ is a linear  \textcolor{black}{space}  and satisfies $\hbox{dim}(\mathcal{D})=3$.
\end{prop}

\noindent{\bf Proof.}
It is obvious that $\mathcal{D}$ is a linear space, and the rest is to prove that the dimension of $\mathcal{D}$ is equal to three. Actually, it follows from Theorem \ref{thm1.1} that for any solution $u=(u_1,u_2,u_3)$ of \eqref{GPs}, there exists a vector $(a_1,a_2,a_3)\in \R^3$ such that $u=(u_1,u_2,u_3)$ satisfies \eqref{exp1}. Applying \eqref{exp1}, we can take a sufficiently large point $-x_0\in\R^+$ such that for $i=1,2,3,$
\begin{equation}\label{non:4.1}
u_i(x_0)=[1+o(1)]a_ie^{\eta_ix_0},\,\ u'_i(x_0)=[1+o(1)]\eta_ia_ie^{\eta_ix_0}\ \ \hbox{as}\,\ -x_0\to +\infty,
\end{equation}
where $\eta_i =\sqrt{|\mu_i|}>0$.

For simplicity, we denote
\begin{equation}\label{26-1}
F_\phi(x_0):=\big(\phi_1(x_0),\phi_2(x_0),\phi_3(x_0)\big),
\end{equation}
and
\begin{equation}\label{26-2}
G_\phi(x_0):=\big(\phi_1(x_0),\phi_2(x_0),\phi_3(x_0),\phi'_1(x_0),\phi'_2(x_0),\phi'_3(x_0)\big).
\end{equation}
We derive from Lemma \ref{lem5.1} that $G_\phi(x_0)$
satisfies the linear  homogeneous  system (\ref{D-1})--(\ref{D-3}), where the coefficient matrix $D$ of $F_\phi(x_0)$
satisfies
 \begin{equation*}\label{non:4.2}
{\footnotesize D=\begin{bmatrix}
[1+o(1)]\mu_1a_1e^{\eta_1x_0}&[1+o(1)]\mu_2a_2e^{\eta_2x_0}&[1+o(1)]\mu_3a_3e^{\eta_3x_0}\\[3mm]
[1+o(1)]\mu_1(\mu_2+\mu_3)a_1e^{\eta_1x_0}&[1+o(1)]\mu_2(\mu_1+\mu_3)a_2e^{\eta_2x_0}&[1+o(1)](\mu_1+\mu_2)\mu_3a_3e^{\eta_3x_0}\\[3mm]
[1+o(1)]\mu_1\mu_2\mu_3a_1e^{\eta_1x_0}&[1+o(1)]\mu_1\mu_2\mu_3a_2e^{\eta_2x_0}&[1+o(1)]\mu_1\mu_2\mu_3a_3e^{\eta_3x_0}\\[3mm]
\end{bmatrix}}
\end{equation*}
as $ -x_0\to +\infty$,  where the exponential decay (\ref{non:4.1}) is used. Since $\mu_1<\mu_2<\mu_3<0$  and $u_i\not\equiv 0$ holds for $i=1,2,3$, we calculate from above that
\[\begin{split}
\det(D)=[1+o(1)]&\mu_1\mu_2\mu_3a_1a_2a_3(\mu_1-\mu_2)(\mu_1-\mu_3)\\
\qquad\qquad &\cdot (\mu_2-\mu_3)e^{(\eta_1+\eta_2+\eta_3)x_0}\neq 0 \quad \mbox{as}\ \,  -x_0\to +\infty,
\end{split}
\]
which further implies that $\hbox{Rank}(D)=3$. Hence, the dimension of the coefficient matrix of $G_\phi(x_0)$
must be equal  to $3$. Therefore, we conclude from above that the dimension of
\[
\begin{aligned}
\big\{&G_\phi(x_0):\,\ G_\phi(x_0)\,\ \hbox{satisfies}\,\ \eqref{Line}\big\}
\end{aligned}
\]
must be equal to $3$. Since $\phi_i(x)$ satisfies the system \eqref{D-1}--\eqref{D-3}, we deduce from (\ref{non:4.1}) that  $\phi_i(x)$ is uniquely determined by $G_\phi(x_0)$. Moreover, we note
that the derivatives of $u_i=u_i(a_1,a_2,a_3,x)$ with respect to $a_i$ provide a three dimensional solutions space, which therefore completes the proof of Proposition \ref{prop5.2}. \qed

Applying Proposition \ref{prop2A}, we are next ready to establish Theorem \ref{thm-nd}.

\vskip 0.1truein

\noindent{\bf Proof of Theorem \ref{thm-nd}.} For convenience, we define $\eta_i=\sqrt{|\mu_i|}>0$ for $i=1, 2, 3$. We address the proof of Theorem \ref{thm-nd} by studying separately the following four different cases:
\begin{equation*}\begin{split}
\mbox{\em Case 1.}\ \   \mu_1<\mu_2<\mu_3<0; \quad
\mbox{\em Case 2.}\ \ \mu_1=\mu_2<\mu_3<0;\\
\mbox{\em Case 3.}\ \ \mu_1<\mu_2=\mu_3<0;\quad
\mbox{\em Case 4.}\ \ \mu_1=\mu_2=\mu_3<0 ;
\end{split}
\end{equation*}
where $u_i\not\equiv 0$ always holds for all $i=1,2,3$.

{\em Case 1.} Suppose $\mu_1<\mu_2<\mu_3<0$. Following  Theorem  \ref{thm1.1},  direct computations yield that there exist at least three linear independent solutions, $i.e.$, three tangent vectors of the solutions manifold with regard to $(a_1,a_2,a_3)$, of the linearized system \eqref{Line} in $H^1(\R)^3$. Thus, the dimension of the solutions for the linearized system \eqref{Line} is at least $3$.  On the other hand, we derive from Lemma \ref{lem5.1} and Proposition \ref{prop5.2} that the dimension of the solutions for the linearized system \eqref{Line} is at most $3$. Therefore, the dimension of the solutions for the linearized system \eqref{Line} must be exactly equal to $3$, and we are done for this case.

{\em Case 2.}  Suppose $\mu_1=\mu_2<\mu_3<0$. We then deduce from Proposition \ref{prop2A} that there exists   $(a_1,a_2,a_3)\in \R^3$ such that $(u_1,u_2,u_3)$ satisfies \eqref{c1} with $a_i\not=0$ for $i=1,2,3$.
In view of \eqref{c1}, we define
\begin{equation}\label{sp-11}
\tilde u_1=u_1\sqrt{1+\big(\frac{u_2}{u_1}\big)^2}=\sqrt{1+\frac{a_2^2}{a_1^2}}\,u_1,
\end{equation}
and
\begin{equation}\label{sp-12}
\tilde\phi_1=\phi_1+\frac{a_2}{a_1}\phi_2,\,\ \tilde\phi_3=\sqrt{1+\frac{a_2^2}{a_1^2}}\,\phi_3.
\end{equation}
One can verify from \eqref{GPs} that $(\tilde u_1,u_3)$ satisfies the following system
\begin{equation}\label{sp-13}
	\left\{\begin{array}{lll}
\tilde u_1''+2(\tilde u_1^{2}+u_3^{2})\tilde u_1=-\mu_1\tilde u_1 \ \ \mbox{in} \ \, \R,\\[3mm]
u_3''+2(\tilde u_1^{2}+u_3^{2})u_3=-\mu_3u_3 \ \ \mbox{in} \ \, \R.
	\end{array}\right.
\end{equation}
Additionally, it follows from \eqref{Line} that $(\tilde \phi_1,\tilde \phi_3)$ satisfies
\begin{equation}\label{sp-14}
	\left\{\begin{array}{lll}
\tilde\phi_1''+2(\tilde u_1^{2}+u_3^{2})\tilde\phi_1+4(\tilde u_1\tilde\phi_1+u_3\tilde\phi_3)\tilde u_1=-\mu_1\tilde\phi_1 \ \ \mbox{in}\, \ \R,\\[3mm]
\tilde\phi_3''+2(\tilde u_1^{2}+u_3^{2})\tilde\phi_3+4(\tilde u_1\tilde\phi_1+u_3\tilde\phi_3)u_3=-\mu_3\tilde\phi_3 \ \ \mbox{in}\, \ \R.
	\end{array}\right.
\end{equation}

Since $\mu_1<\mu_3<0$ and $(\tilde u_1,u_3)\neq (0,0)$, the same argument as above (i) implies that $(\tilde u_1,u_3)$ is a non-degenerate critical point of the system \eqref{sp-13}, and the dimension of solutions $(\tilde \phi_1,\tilde \phi_3)$ for \eqref{sp-14} is two. Let $v_1=(v_{11},v_{13})$ and $v_2=(v_{21},v_{23})$ be the two linearly independent solutions of  \eqref{sp-14}. Thus, $(\tilde \phi_1,\tilde \phi_3)$ solves \eqref{sp-14}, if and only if
\begin{equation}\label{sp-15}
(\tilde \phi_1,\tilde \phi_3)=a(v_{11},v_{13})+b(v_{21},v_{23})\ \ \mbox{holds for some}\, \ a,\, b\in\R.
\end{equation}
Moreover, calculating $\eqref{Line}_2-\frac{a_2}{a_1}\times\eqref{Line}_1$, it yields that
\begin{equation}\label{sp-18}
\Big(\phi_2-\frac{a_2}{a_1}\phi_1\Big)''+2(u_1^{2}+u_2^{2}+u_3^{2})\Big(\phi_2-\frac{a_2}{a_1}\phi_1\Big)
=-\mu_1\Big(\phi_2-\frac{a_2}{a_1}\phi_1\Big) \ \ \mbox{in} \ \, \R.
\end{equation}
We then follow from
\eqref{sp-11}, \eqref{sp-12} and \eqref{sp-15} that \eqref{Line} admits the following two linearly independent solutions
\begin{equation}\label{sp-16} (\phi_1,\phi_2,\phi_3)=\Big(\frac{a_1^2}{a_1^2+a_2^2}v_{11},\frac{a_1a_2}{a_1^2+a_2^2}v_{11},
\frac{a_1}{\sqrt{a_1^2+a_2^2}}v_{13}\Big)
\end{equation}
and
\begin{equation}\label{sp-17} (\phi_1,\phi_2,\phi_3)=\Big(\frac{a_1^2}{a_1^2+a_2^2}v_{21},\frac{a_1a_2}{a_1^2+a_2^2}
v_{21},\frac{a_1}{\sqrt{a_1^2+a_2^2}}v_{23}\Big).
\end{equation}

To seek for other nontrivial solutions $(\phi_1,\phi_2,\phi_3)$ of \eqref{Line}, we now assume that
\begin{equation}\label{7:sp-17}
(\tilde \phi_1,\tilde \phi_3)=\Big(\phi_1+\frac{a_2}{a_1}\phi_2,\sqrt{\frac{a_1^2+a_2^2}{a_1^2}}\phi_3\Big)=(0,0),
\end{equation}
which also implies that $\phi_3\equiv 0$. Since it follows from \eqref{c1} that $u_1$ cannot change sign, we deduce that $\mu_1<0$ is the first eigenvalue of $L_u$ defined in (\ref{1:op}). It then yields from (\ref{1:op}) and \eqref{sp-18} that
$\phi_2-\frac{a_2}{a_1}\phi_1\equiv Cu_1$ holds for some $C_1\in\R$. If $C_1=0$, $i.e.$, $\phi_2-\frac{a_2}{a_1}\phi_1\equiv 0$, then we derive from (\ref{7:sp-17}) that $(\phi_1,\phi_2,\phi_3)=(0,0,0)$, which is however trivial. Thus,  $\phi_2-\frac{a_2}{a_1}\phi_1\equiv Cu_1$ holds for some $C_1\neq 0$. We then deduce from (\ref{7:sp-17}) that the third solution $(\phi_1,\phi_2,\phi_3)$ of \eqref{Line} must satisfy
\begin{equation}\label{sp-19}
(\phi_1,\phi_2,\phi_3)=k\Big(-\frac{a_2}{a_1}u_1,u_1,0\Big)=k(-u_2,u_1,0)\,\ \hbox{for some}\,\ k\neq 0.
\end{equation}
One can verify that those three solutions   (\ref{sp-16}), (\ref{sp-17}) and (\ref{sp-19})  of \eqref{Line} are linearly independent. Further, it follows from (\ref{sp-14}), (\ref{sp-15}) and (\ref{sp-19}) that any solution of \eqref{Line} can be expressed as a linear combination of these three independent solutions.
Hence, the dimension of solutions $(\phi_1,\phi_2,\phi_3)$ for \eqref{Line} is  exactly three.  This  proves the non-degeneracy of the system \eqref{GPs} in this case.

{\em Case 3.}  Suppose  $\mu_1<\mu_2=\mu_3<0$. \textcolor{black}{The proof of Case 3 is very similar to that of Case 2, except that the  eigenvalue $\mu_2$ of  $L_u$ in $L^2(\R)$ may be multiple, which would lead to the degeneracy of solutions for the system \eqref{Lines}. In the following, hence we only need to prove that $\mu_2$ is also a simple eigenvalue of $L_u$.}

In order to prove $\mu_2$ is a simple eigenvalue of $L_u$,  we first claim that
for any $0<\varepsilon <1$, there exists a sufficiently large constant $\bar N>0$  such that
\begin{equation}\label{ob-8}
\Big|\frac{u_0(x_1)}{e^{\eta_2x_1}}-\frac{u_0(x_2)}{e^{\eta_2x_2}}\Big|<\varepsilon\ \ \hbox{for any}\ \ x_2<x_1<-\bar N,
\end{equation}
where $u_0(x)\not\equiv 0$ is any solution of $L_u=\mu_2 u$ in $L^2(\R)$, and the operator $L_u$ is defined in (\ref{1:op}).

To prove (\ref{ob-8}), we note from Theorem \ref{thm1.1} that
$$2(u_1^2+u_2^2+u_3^2)=2[1+o(1)]a_1^2e^{2\eta_1x}\ \ \mbox{ as}\ \ x\to-\infty,$$
which implies that there exists a sufficiently large constant $N_1>0$ such that
\begin{equation}\label{ob-2}
2(u_1^2+u_2^2+u_3^2)\leq 4a_1^2e^{2\eta_1x},\,\ \forall x<-N_1.
\end{equation}
Since $u_0(x)\not\equiv 0$ is a solution of $L_u=\mu_2 u$ in $L^2(\R)$, we can choose a sufficiently large value $-x_0>N_1$ such that $u_0(x_0)\neq 0$. Moreover,  for simplicity we may assume $u_0(x_0)>0$. Taking $U_1(x):=u_0(x_0)e^{-\eta_2x_0}e^{\eta_2x}$ as a comparison function, we have $U_1(x_0)=u_0(x_0)$ and
\[
\begin{aligned}
(L_u-\mu_2)U_1&=-\frac{d^2}{dx^2}U_1-\mu_2U_1-2(u_1^2+u_2^2+u_3^2)U_1\\
&=-2(u_1^2+u_2^2+u_3^2)U_1<0\ \ \hbox{in}\,\ \R.	
\end{aligned}
\]
By the comparison principle, we then obtain from above that
\begin{equation}\label{ob-3}
u_0(x)\geq U_1(x)=u_0(x_0)e^{-\eta_2x_0}e^{\eta_2x},\,\   x\leq x_0,
\end{equation}
where $-x_0 >N_1 $ is sufficiently large.

On the other hand, we  define for any $0<\varepsilon <1$,
\[
U_2(x):=\Big[u_0(x_0)e^{-\eta_2x_0}+\frac{\varepsilon}{2}\Big]e^{\eta_2x}-Be^{(\eta_2+2\eta_1)x},
\]
where  $B=a_1^2\eta_1^{-2}\big(\|u_0\|_{L^\infty(\R)}e^{-\eta_2x_0}+1\big)>0$, and the sufficiently large constant $-x_0 >N_1 $ is  as in (\ref{ob-3}). We then have for any $0<\varepsilon <1$,
\[
\begin{aligned}
U_2(x_0)&=u_0(x_0)+\frac{\varepsilon}{2}e^{\eta_2x_0}-a_1^2\eta_1^{-2}\Big[\|u_0\|_{L^\infty(\R)}e^{-\eta_2x_0}+
\frac{\varepsilon}{2}\Big]e^{(\eta_2+2\eta_1)x_0}\\
&\geq u_0(x_0)+\frac{\varepsilon}{2}e^{\eta_2x_0}-a_1^2\eta_1^{-2}\Big[\|u_0\|_{L^\infty(\R)}e^{-\eta_2x_0}+1\Big]e^{(\eta_2+2\eta_1)x_0}\\
&=u_0(x_0)+e^{\eta_2x_0}\Big\{\frac{\varepsilon}{2}-a_1^2\eta_1^{-2}\Big[\|u_0\|_{L^\infty(\R)}e^{(2\eta_1-\eta_2)x_0}
+e^{2\eta_1x_0}\Big]\Big\}.
\end{aligned}
\]
Since $2\eta_1-\eta_2>0$ and $\eta_1>0$, there exists a sufficiently large constant $N_2>0$ such that
\[
\frac{\varepsilon}{2}-a_1^2\eta_1^{-2}\Big[\|u_0\|_{L^\infty(\R)}e^{(2\eta_1-\eta_2)x_0}+e^{2\eta_1x_0}\Big]>0\,\ \hbox{for all}\,\ x_0<-N_2,
\]
which implies that $U_2(x_0)>u_0(x_0)$ holds for all $x_0<-N_2$. Moreover, we obtain from \eqref{ob-2} that
\[
\begin{aligned}
	&(L_u-\mu_2)U_2\\
 >&-4a_1^2\big[u_0(x_0)e^{-\eta_2x_0}+\frac{\varepsilon}{2}\big]
 e^{(\eta_2+2\eta_1)x}+B\big[(\eta_2+2\eta_1)^2-\eta_2^2\big]e^{(\eta_2+2\eta_1)x}\\
	>&4\big[B\eta_1^2-a_1^2\big(\|u_0\|_{L^\infty(\R)}e^{-\eta_2x_0}+1\big)\big]e^{(\eta_2+2\eta_1)x}=0\ \ \mbox{in}\ \, \R.
\end{aligned}
\]
By the comparison principle, we then deduce from above that for all $0<\varepsilon <1$,
\begin{equation}\label{ob-4}
\begin{split}
u_0(x)\leq U_2(x)=&\Big[u_0(x_0)e^{-\eta_2x_0}+\frac{\varepsilon}{2}\Big]e^{\eta_2x}\\
&-a_1^2\eta_1^{-2}\big(\|u_0\|_{L^\infty(\R)}e^{-\eta_2x_0}+1\big)e^{(\eta_2+2\eta_1)x}, \ \ \forall x\leq x_0,
\end{split}
\end{equation}
where $-x_0>\max\big\{N_1,N_2\big\}$ is sufficiently large.
We now conclude from \eqref{ob-3} and \eqref{ob-4} that  for all $0<\varepsilon <1$,
\begin{equation}\label{ob-5}
\begin{split}
u_0(x_0)e^{-\eta_2x_0}\leq \frac{u_0(x)}{e^{\eta_2x}}&\leq u_0(x_0)e^{-\eta_2x_0}+\frac{\varepsilon}{2}-a_1^{-2}\big(\|u_0\|_{L^\infty(\R)}e^{-\eta_2x_0}+1\big)e^{2\eta_1x}\\
&\le u_0(x_0)e^{-\eta_2x_0}+\frac{\varepsilon}{2},\ \ \mbox{if}\ \, x\leq x_0<-\max\big\{N_1,N_2\big\},
\end{split}
\end{equation}
which therefore gives that \eqref{ob-8} hold true.

By Cauchy's convergence criteria, we further derive from  \eqref{ob-8} that if $u_0(x)\not\equiv 0$ is a nonzero solution of $L_u=\mu_2 u$ in $L^2(\R)$, where the operator $L_u$ is defined in (\ref{1:op}), then $u_0(x)$ satisfies
\begin{equation}\label{ob-1}
\lim_{x\to-\infty}\frac{u_0(x)}{e^{\eta_2x}}=A\in \R\backslash\{0\}.
\end{equation}
Suppose now that $v_1\not\equiv 0$ and $v_2\not\equiv 0$ are two different solutions of $(L_u-\mu_2)v=0$ in $\R$. Following \eqref{ob-1}, we may assume that
\[
\lim_{x\to-\infty}\frac{v_i(x)}{e^{\eta_2x}}=A_i\neq 0,\,\ \hbox{where}\,\ i=1,2.
\]
Set $v_0=v_1-\frac{A_1}{A_2}v_2$, which yields from above that $\lim_{x\to-\infty}\frac{v_0(x)}{e^{\eta_2x}}=0$ holds for $\eta_2>0$. This further implies that $v_0\equiv 0$ in $\R$. We therefore conclude that $\mu_2$ is a simple eigenvalue of $L_u$ in $\R$, where the operator $L_u$ is defined in (\ref{1:op}).

\textcolor{black}{Since \eqref{ob-8} holds true, we can proceed as Case 2 to further obtain the non-degeneracy of Case 3, and we are therefore done.}

{\em Case 4.}  Suppose $\mu_1=\mu_2=\mu_3<0$. Following Proposition \ref{prop2A} (iii), $u_i$ never change sign, and $u_i$, $u_j$ are proportional, where $i,j=1,2,3$. Define
\begin{equation}\label{sp-26}
	\tilde u=u_1\sqrt{1+\big(\frac{u_2}{u_1}\big)^2+\big(\frac{u_3}{u_1}\big)^2}\,u_1,
\end{equation}
and
\begin{equation}\label{sp-27}
	\tilde\phi=\phi_1+\frac{u_2}{u_1}\phi_2+\frac{u_3}{u_1}\phi_3,
\end{equation}
then Proposition \ref{prop2A} (iii) yields that
\begin{equation}\label{sp-28}
	\tilde u''+2\tilde u^3=-\mu_1\tilde u\ \ \mbox{in}\ \, \R,
\end{equation}
and
\begin{equation}\label{sp-29}
	\tilde \phi''+4\tilde u^2\tilde \phi=-\mu_1\tilde \phi\ \ \mbox{in}\ \, \R.
\end{equation}
The nondegeneracy of \eqref{sp-29} yields that there exists a constant $C_0\in\R$ such that
\begin{equation}\label{O3}
\tilde \phi=\phi_1+\frac{u_2}{u_1}\phi_2+\frac{u_3}{u_1}\phi_3=C_0\tilde u'.
\end{equation}
Moreover, similar to \eqref{sp-18}, we have
\begin{equation}\label{O4}
\phi_j-\frac{a_j}{a_1}\phi_1=C_ju_1, \ \ j=1,2.	
\end{equation}

Following \eqref{O3} and \eqref{O4}, we conclude that the system \eqref{Line} admits the following three linearly independent solutions:
\[
(\phi_1,\phi_2,\phi_3)=(-u_2,u_1,0),\,\ (-u_2,0,u_3)\,\ \hbox{or}\,\ (u'_1,u'_2,u'_3),
\]
and any solution of \eqref{Line} can be expressed as a linear combination of the above three solutions.
This further proves the non-degeneracy of Case 4, which therefore completes the proof of Theorem \ref{thm-nd}.  \qed

\section{The General Case of $N>3$}
In this section, we discuss whether all main results of the present paper can be extended generally to the  $N-$component Schr\"{o}dinger system (\ref{GPN}), where $\mu_1\leq\mu_2\leq\cdots\leq\mu_N<0$ and $N>3$. We first have the following similar version of Lemma \ref{lem2.2}.

\begin{lem}\label{lem5-1}
Suppose $u=(u_1,u_2,\cdots,u_N)\in H^1(\R)^N$ is a   solution of \eqref{GPN}, where $\mu_1\leq\mu_2\leq\cdots\leq\mu_N<0$ and $N>3$. Then we have the following $N$ constants of motion:	
\begin{equation}\label{N.1}
\begin{split}
			\sum_{i=1}^N(u_i')^2+ \Big(\sum_{i=1}^Nu_i^2\Big)^2+\sum_{i=1}^N\mu_iu_i^2=0 \ \ \mbox{in}\, \ \R,	
\end{split}
\end{equation}
\begin{equation}\label{N.2}
\begin{split}
&\sum_{1\leq j_1<j_2\leq N}(u_{j_1}'u_{j_2}-u_{j_1}u_{j_2}')^2+\Big(\sum_{i=1}^Nu_i^2\Big)\sum_{1\leq j_1,\,j_2\leq N,\,j_1\neq j_2}\mu_{j_1}u^2_{j_2}\\[2mm]
&\,\,+\sum_{1\leq j_1,\,j_2\leq N,\,j_1\neq j_2} \mu_{j_1}(u'_{j_2})^2+\sum_{1\leq j_1,\,j_2\leq N,\,j_1\neq j_2} \mu_{j_1}\mu_{j_2}u^2_{j_2}=0 \ \ \mbox{in}\, \ \R,
		\end{split}
\end{equation}
and the $k$-th ($3\leq k\leq N$) identity satisfies
\begin{equation}\label{N.3}
\begin{split}
&\sum_{\stackrel{ 1\leq j_1<j_2<\cdots<j_{k-2}\leq N,\, 1\leq j_{ k-1}<j_{k}\leq N,}
{ \prod_{1\leq l<m\leq k}(j_l-j_m)\neq 0  } }\mu_{j_1}\mu_{j_2}\cdots\mu_{j_{k-2}}(u_{j_{k-1}}'u_{j_k}-u_{j_{k-1}}u_{j_k}')^2\\[2mm]
			&\,\,+\Big(\sum_{i=1}^Nu_i^2\Big)\sum_{1\leq j_1<j_2<\cdots<j_{k-1}\leq N,\,1\leq j_{k}\leq N,\,\prod_{1\leq l<m\leq k}(j_l-j_m)\neq 0 }\mu_{j_1}\mu_{j_2}\cdots\mu_{j_{k-1}}u_{j_k}^2\\[2mm]
			&\,\,+\sum_{1\leq j_1<j_2<\cdots<j_{k-1}\leq N,\,1\leq j_{k}\leq N,\,\prod_{1\leq l<m\leq k}(j_l-j_m)\neq 0 } \mu_{j_1}\mu_{j_2}\cdots\mu_{j_{k-1}}(u'_{j_k})^2\\[2mm]
			&\,\,+\sum_{1\leq j_1<j_2<\cdots<j_{k-1}\leq N,\,1\leq j_{k}\leq N,\,\prod_{1\leq l<m\leq k}(j_l-j_m)\neq 0 } \mu_{j_1}\mu_{j_2}\cdots\mu_{j_{k-1}}\mu_{j_{k}}u^2_{j_k}=0 \ \ \mbox{in}\, \ \R.
		\end{split}
	\end{equation}
\end{lem}

\noindent{\bf Proof.} The process of proving Lemma \ref{lem5-1} is very similar to that of Lemma \ref{lem2.2}. Especially, the $k$-th  ($3\leq k\le N$)  identity (\ref{N.3}) is obtained through multiplying the $j$-th of \eqref{GPs} by the following integral factor
\[
\begin{aligned}
	\sum_{\stackrel{ 1\leq j_1<j_2<\cdots<j_{k-2}\leq N,\,1\leq j_{k-1}\leq N,}{\prod_{1\leq l<m\leq k}(j_l-j_m)(j_l-j)(j_m-j)\neq 0}}&2\Big[\mu_{j_1}\mu_{j_2}\cdots\mu_{j_{k-2}}(u'_{j}u_{j_{k-1}}-u'_{j_{k-1}}u_{j})u_{j_{k-1}}\\
	&\quad+\frac{\mu_{j_1}\mu_{j_2}\cdots\mu_{j_{k-1}}}{k-1}u'_{j}\Big].
\end{aligned}
\]
We leave the detailed proof to the interested reader. \qed

In view of Lemma \ref{lem5-1}, we then have the following two conjectures, which are expected to be involved with more complicated algebraic analysis.

\vskip 0.05truein
\noindent
{\em {\bf Conjecture 5.1.}
Suppose $u=(u_1,u_2,\cdots,u_N)\in H^1(\R)^N$ is a   solution of \eqref{GPN}, where $\mu_1\leq\mu_2\leq\cdots\leq\mu_N<0$ and $N>3$. Then there exists a unique vector $(a_1,a_2,\cdots,a_N)\in\R^N$ such that the solution $u$  satisfies
\begin{equation}\label{5:N.3}
	u_i(x)=\frac{g_{i}(x)}{f(x)},\quad i=1,2,\cdots,N,
\end{equation}
where $g_i(x)$ and $f(x)$ satisfy for $\eta_i=\sqrt{|\mu_i|}>0$,
\[
\begin{aligned}
		g_i(x)&=a_ie^{\eta_ix}+\sum_{j=1}^{N}\frac{a_ia_j^2(\eta_i-\eta_j)}{4\eta_j^2(\eta_i+\eta_j)}e^{(\eta_i+2\eta_j)x}\\
		&\quad+\sum_{1\leq j_1<j_2\leq N}\frac{a_ia_{j_1}^2a_{j_2}^2(\eta_{j_1}-\eta_{j_2})^2(\eta_i-\eta_{j_1})(\eta_i-\eta_{j_2})}{4^2\eta_{j_1}^2\eta_{j_2}^2(\eta_{j_1}+\eta_{j_2})^2(\eta_i+\eta_{j_1})(\eta_i+\eta_{j_2})}e^{(\eta_i+2\eta_{j_1}+2\eta_{j_2})x}\\
		&\quad+\cdots\cdots\quad\\
		&\quad+\sum_{1\leq j_1<j_2<\cdots<j_k\leq N}\frac{1}{4^{k}}\prod_{l=1}^k\frac{a_ia_{j_l}^2(\eta_i-\eta_{j_l})}{\eta_{j_l}^2(\eta_i+\eta_{j_l})}\prod_{1\leq l<m\leq k}\frac{(\eta_{j_l}-\eta_{j_m})^2}{(\eta_{j_l}+\eta_{j_m})^2}e^{(\eta_i+2\eta_{j_1}+2\eta_{j_2}+\cdots+2\eta_{j_k})x}\\
		&\quad+\cdots\cdots\\
		&\quad+\sum_{1\leq j_1<j_2<\cdots<j_{N-1}\leq N}\frac{1}{4^{N-1}}\prod_{l=1}^N\frac{a_ia_{j_l}^2(\eta_i-\eta_{j_l})}{\eta_{j_l}^2(\eta_i+\eta_{j_l})}\\
&\qquad\qquad \cdot\prod_{1\leq l<m\leq N}\frac{(\eta_{j_l}-\eta_{j_m})^2}{(\eta_{j_l}+\eta_{j_m})^2}e^{(\eta_i+2\eta_{j_1}+2\eta_{j_2}+\cdots+2\eta_{j_N})x}, \ \ i=1, 2, \cdots, N,
\end{aligned}
\]
and
\[
\begin{aligned}
		f(x)&=1+\sum_{j=1}^{N}\frac{a_j^2}{4\eta^2_j}e^{2\eta_jx}+\sum_{1\leq j_1<j_2\leq N}\frac{a_{j_1}^2a_{j_2}^2(\eta_{j_1}-\eta_{j_2})^2}{4^2\eta_{j_1}^2\eta_{j_2}^2(\eta_{j_1}+\eta_{j_2})^2}e^{(2\eta_{j_1}+2\eta_{j_2})x}\\
		&\quad+\cdots\cdots\quad\\
		&\quad+\sum_{1\leq j_1<j_2<\cdots<j_k\leq N}\frac{1}{4^k}\prod_{1\leq l<m\leq k}\frac{a_{j_l}^2(\eta_{j_l}-\eta_{j_m})^2}{\eta_{j_l}^2(\eta_{j_l}+\eta_{j_m})^2}e^{(2\eta_{j_1}+2\eta_{j_2}+\cdots+2\eta_{j_k})x}\\
		&\quad+\cdots\cdots\\
		&\quad+\sum_{j_1=1,\,j_2=2,\,\cdots,\,j_N=N}\frac{1}{4^N}\prod_{1\leq l<m\leq N}\frac{a_{j_l}^2(\eta_{j_l}-\eta_{j_m})^2}{\eta_{j_l}^2(\eta_{j_l}+\eta_{j_m})^2}e^{( 2\eta_{j_1}+2\eta_{j_2}+\cdots+2\eta_{N})x}\, _.
	\end{aligned}
	\]
}

\noindent
{\em {\bf Conjecture 5.2.} For $\mu_1\leq\mu_2\leq\cdots\leq\mu_N<0$, where $N>3$, \eqref{GPN} does not admit any orthonormal solution $u=(u_1,u_2,\cdots,u_N)\in H^1(\R)^N$ satisfying $(u_i, u_j)=\delta_{ij}$ for all $1\leq i,j\leq N$.
}

We remark that Conjecture 5.2 already appeared earlier in \cite[Section 1.3]{Lewin}.
For any solution $u=(u_1,u_2,\cdots,u_N)$ of \eqref{GPN}, we next consider the following linearized system of \eqref{GPN} around the solution $u$:
\begin{equation}\label{GPNLN}
	\begin{array}{lll}
		u_i''+2\Big(\sum_{k=1}^Nu_k^2\Big)\phi_i+4\Big(\sum_{k=1}^Nu_k\phi_k\Big)u_i=-\mu_i\phi_i \ \,\ \mbox{in}\, \ \R , \ \  i=1, 2, \cdots, N.
	\end{array}
\end{equation}
If Conjecture 5.1 holds true, then we are able to prove the following non-degeneracy of (\ref{GPNLN}).

\begin{thm}\label{thm5.2} Suppose $u=(u_1,u_2,\cdots,u_N)\in H^1(\R)^N$ is a    solution of \eqref{GPN} satisfying $u_i\not\equiv 0$ for $i=1,2,\cdots,N$, where $\mu_1<\mu_2<\cdots<\mu_N<0$ and $N>3$. If Conjecture 5.1 holds true, then $u$ is non-degenerate, in the sense that the dimension of solutions  for the linearized system \eqref{GPNLN} in $H^1(\R)^N$ is equal to $N$.
\end{thm}

\noindent{\bf Proof.}  Similar to Lemma \ref{lem5.1}, one can derive that $(\phi_1,\phi_2,\cdots,\phi_N)$ satisfies the following $N$ constants of motion, which correspond to the linearized version of Lemma \ref{lem5-1}:
	\begin{equation}\label{L.1}
	\begin{split}
		\sum_{i=1}^Nu_i'\phi_i'+ 2\Big(\sum_{i=1}^Nu_i^2\Big)\Big(\sum_{i=1}^Nu_i\phi_i\Big)+\sum_{i=1}^N\mu_iu_i\phi_i=0 \ \,\ \mbox{in}\, \ \R,	
	\end{split}
\end{equation}
\begin{equation}\label{L.2}
	\begin{split}
		&\sum_{1\leq j_1<j_2\leq N}\big(u_{j_1}'u_{j_2}-u_{j_1}u_{j_2}'\big)\big(\phi_{j_1}'u_{j_2}+u_{j_1}'\phi_{j_2}-\phi_{j_1}u_{j_2}'-u_{j_1}\phi_{j_2}'\big)\\[2mm]
		+&\sum_{1\leq j_1,\,j_2\leq N,\,j_1\neq j_2}\Big[\Big(\sum_{i=1}^Nu_i^2\Big)\mu_{j_1}u_{j_2}\phi_{j_2}+\Big(\sum_{i=1}^Nu_i\phi_i\Big)\mu_{j_1}u_{j_2}^2\Big]\\[2mm]
		+&\sum_{1\leq j_1,\,j_2\leq N,\,j_1\neq j_2} \mu_{j_1}u'_{j_2}\phi'_{j_2}+\sum_{1\leq j_1,j_2\leq N,j_1\neq j_2} \mu_{j_1}\mu_{j_2}u_{j_2}\phi_{j_2}=0 \ \,\ \mbox{in}\, \ \R,
	\end{split}
\end{equation}
and the $k$-th ($3\leq k\leq N$) identity satisfies
\begin{equation}\label{L.3}
	\begin{split}
		&\sum_{\stackrel{ 1\leq j_1<j_2<\cdots<j_{k-2}\leq N,\,1\leq j_{k-1}<j_{k}\leq N,}{\prod_{1\leq l<m\leq k}(j_l-j_m)\neq 0 } }\mu_{j_1}\mu_{j_2}\cdots\mu_{j_{k-2}}(u_{j_{k-1}}'u_{j_k}-u_{j_{k-1}}u_{j_k}')\\[2mm]		 &\qquad\qquad\qquad\quad\cdot\Big[(\phi_{j_{k-1}}'u_{j_k}+u_{j_{k-1}}'\phi_{j_k}-\phi_{j_{k-1}}u_{j_k}'-u_{j_{k-1}}\phi_{j_k}')\Big]\\[2mm]
		+&\sum_{i=1}^Nu_i^2\sum_{\stackrel{ 1\leq j_1<j_2<\cdots<j_{k-1}\leq N,\,1\leq j_{k}\leq N,}{\prod_{1\leq l<m\leq k}(j_l-j_m)\neq 0} }\mu_{j_1}\mu_{j_2}\cdots\mu_{j_{k-1}}u_{j_k}\phi_{j_k}\\[2mm]
		+&\sum_{i=1}^Nu_i\phi_i\sum_{\stackrel{ 1\leq j_1<j_2<\cdots<j_{k-1}\leq N,\,1\leq j_{k}\leq N,}{\prod_{1\leq l<m\leq k}(j_l-j_m)\neq 0} }\mu_{j_1}\mu_{j_2}\cdots\mu_{j_{k-1}}u_{j_k}^2\\[2mm]
		+&\sum_{\stackrel{ 1\leq j_1<j_2<\cdots<j_{k-1}\leq N,\,1\leq j_{k}\leq N,}{\prod_{1\leq l<m\leq k}(j_l-j_m)\neq 0 } }\mu_{j_1}\mu_{j_2}\cdots\mu_{j_{k-1}}(u'_{j_k}\phi'_{j_k})\\[2mm]
		+&\sum_{\stackrel{ 1\leq j_1<j_2<\cdots<j_{k-1}\leq N,\,1\leq j_{k}\leq N,}{\prod_{1\leq l<m\leq k}(j_l-j_m)\neq 0 }} \mu_{j_1}\mu_{j_2}\cdots\mu_{j_{k-1}}\mu_{j_{k}}u_{j_k}\phi_{j_k}=0 \ \,\ \mbox{in}\, \ \R.
	\end{split}
\end{equation}

If Conjecture 5.1 holds true, then we can take a sufficiently large point $-x_0\in\R^+$ so that
\begin{equation}
	u_i(x_0)=[1+o(1)]a_ie^{\eta_ix_0},\,\ u'_i(x_0)=[1+o(1)]\eta_ia_ie^{\eta_ix_0}\ \ \hbox{as}\,\ -x_0\to +\infty, \ \
\end{equation}
where $i=1,2,\cdots,N$.
Applying \eqref{L.1}--\eqref{L.3}, the same argument of Proposition \ref{prop5.2}  then yields that Theorem \ref{thm5.2} holds true, and we are done. \qed

We finally comment that the proof of Theorem \ref{thm5.2} depends strongly on the assumption $\mu_1<\mu_2<\cdots<\mu_N<0$. Stimulated by  Theorem \ref{thm-nd}, we guess that if $\mu_j=\mu_i$ holds for some $i\neq j$, then the non-degeneracy of Theorem \ref{thm5.2} is still true for all $N>3$, but its proof is more involved with the inductive analysis, which is also left to the interested reader.

\vskip 0.5truein

\noindent {\bf Acknowledgements:} The authors are very grateful to the reviewers for their valuable suggestions and comments, which lead to the great improvements of the present paper. The authors also thank Professors R. Frank, D. Gontier and M. Lewin very much for their helpful discussions on the present paper.

\end{document}